\documentclass[11pt]{llncs}
\usepackage[letterpaper,hmargin=0.5in,vmargin=1.15in]{geometry}
\usepackage{amsmath,amssymb,amsbsy,amsfonts,latexsym,amsopn,amstext,amsxtra,euscript,amscd,color}
\usepackage{array}
\usepackage{booktabs}
\usepackage{graphicx}
\usepackage{slashbox}
\usepackage{psfrag}
\usepackage{epsfig}
\usepackage[sort&compress,numbers,square]{}
\def\00{{\bf 0}}
\def\+{\oplus}

\def\\{\cr}
\def\({\left(}
\def\){\right)}

\providecommand{\newoperator}[3]{%
  \newcommand*{#1}{\mathop{#2}#3}}

\newoperator{\FD}{\mathrm{FD}}{\nolimits}

\begin{document}
\title{A new numerical algorithm for the solutions of hyperbolic partial
differential equations in $(2+1)$-dimensional space}
\author{Brajesh Kumar Singh}
\institute{Department of Applied Mathematics, School for Physical Sciences, \\ Babasaheb Bhimrao Ambedkar University Lucknow-226 025 (UP)  INDIA \\
\email{bksingh0584@gmail.com}}
\date{\today}

\maketitle
\begin{abstract}
This paper deals with a construction of new algorithm: the modified trigonometric cubic B-Spline differential quadrature (MTB-DQM) for space discretization together with a time integration algorithm" for numerical computation of the hyperbolic equations. Specially, MTB-DQM has been implemented for the initial value system of the telegraph equations together with both Dirichlet and Neumann type boundary conditions. The MTB-DQM is a DQM based on modified trigonometric cubic B-splines as new  base functions. The problem has been reduced into an amenable system of ordinary differential equations adopting MTB-DQM. The resulting system of ordinary differential equations is solved using time integration algorithms. Further, the stability of MTB-DQM is studied by computing the eigenvalues of the coefficients matrices for various grid points, which confirmed the stability of MTB-DQM for the telegraphic equations. The accuracy of the method has been illustrated in terms of the various discrete error norms for six test problems of the telegraph equation. A comparison of computed numerical solutions with that obtained by the other methods has been carried out for various time levels considering various space sizes.
\end{abstract}

%{\bf In this paper, an algorithm of differential quadrature method based on modified trigonometric cubic B-Spline base functions in space has been developed to transform a partial differential equation into time dependent ordinary differential equation. The well known     }


\section{Introduction}
The hyperbolic partial differential equations have a great attention due to remarkable applications in fields of applied science and engineering, for instance, these equations model fundamental equations in atomic physics \cite{LS10} and are very useful in understanding various physical phenomena in applied sciences and engineering. It models the vibrations of structures (e.g. buildings, machines and beams). Consider the initial value linear system of second-order hyperbolic telegraph equation in $(2+1)$ dimension together with the following boundary conditions:
\begin{equation}\label{int-tel-eqn1}
\left\{\begin{split}
 & \frac{\partial^2 u (x, y, t)}{\partial t^2} +2\alpha  \frac{\partial u (x, y, t)}{\partial t} +\beta^2 u(x,y,t)=  \frac{\partial^2 u (x, y, t)}{\partial x^2} + \frac{\partial^2 u (x, y, t)}{\partial y^2}+f(x,y,t),~~ \\
 & u(x,y,0)=\phi(x,y), \qquad  u_t(x,y,0)=\psi(x,y), ~~~~ (x,y)\in \Omega, t>0.
 \end{split}\right.
\end{equation}
$a)$ Dirichlet boundary condition:
\begin{equation}\label{eqn-DBC}
 \begin{array}{ll}
 u(0,y,t)=\phi_1(y,t),u(1,y,t)=\phi_2(y,t),
 u(x,0,t)=\phi_3(x,t),u(x,1,t)=\phi_4(x,t),
\end{array} (x,y)\in \partial \Omega, t>0
\end{equation}
\noindent or  $b)$ Neumann boundary conditions:
\begin{equation}\label{eqn-NBC}
\begin{array}{ll}
\frac{\partial u}{\partial x}(0,y,t)=\psi_1(y,t), \frac{\partial u}{\partial x}(1,y,t)=\psi_2(y,t), \frac{\partial u}{\partial y}(x,0,t)=\psi_3(x,t), \frac{\partial u}{\partial y}(x,1,t)=\psi_4(x,t),
\end{array} (x,y)\in \partial \Omega, t>0.
\end{equation}
where $\partial \Omega$ denotes the boundary of the computational domain $\Omega=[0,1]\times [0,1] \subset R^2$ and $\alpha>0, \beta$ are arbitrary constants, and $\psi, \phi, \psi_i, \phi_i (i=1,2 ,3,4)$ are known smooth functions. Eq. \eqref{int-tel-eqn1} with $\beta=0$ is a damped wave equation while for $\beta >0$ it is known as hyperbolic telegraph equation. This equation is more convenient than ordinary diffusion equation in modeling reaction diffusion for such branches of sciences \cite{DG10}, and mostly used in wave propagation of electric signals in a cable transmission line \cite{P90}.

In the recent years, various type of numerical techniques have been developed for solving partial differential equations \cite{ED15,AMS14,SS14,SA14,SAS16,SB16,AS13,SPK16}, among others, initial value problems of telegraph equation, in one dimension, have been solved by Taylor matrix method \cite{BS11}, dual reciprocity boundary integral method \cite{DG10}, unconditionally stable finite difference scheme \cite{GC07}, implicit difference scheme \cite{M09}, variational iteration method \cite{MYL11}, modified B-spline collocation method \cite{MB13}, Chebyshev tau method \cite{SD10}, interpolating scaling function method \cite{LS10},  cubic B-spline collocation method \cite{SR16} whereas two dimensional initial value problem of the telegraph equations have been solved by Taylor matrix method by B\"{u}lb\"{u}l and Sezer \cite{BS11a} which  converts the telegraph equation into the matrix equation, Two meshless methods- namely meshless local weak-strong (MLWS) and meshless local Petrov-Galerkin (MLPG) method by Dehghan and Ghesmati \cite{DG10a}, higher order implicit collocation method \cite{DM09}, A polynomial based differential quadrature method \cite{JPM12}, modified cubic B-spline differential quadrature method \cite{MB14}, an unconditionally stable alternating direction implicit scheme \cite{MJ01}, A hybrid method due to Dehghan and Salehi \cite{DS12}, compact finite difference scheme by Ding and Zhang \cite{DZ09} with accuracy of order four in both space and time. Two dimensional linear hyperbolic telegraph equation with variable coefficients has been solved by Dehghan and Shorki \cite{DS09}.

The differential  quadrature method (DQM) has been developed by Bellman et al. \cite{BKLV75,BKC72} to obtain the numerical solutions of partial differential equations (PDEs), have received a great attention among the researchers. After seminal work of Bellman et al.\cite{BKLV75,BKC72}, and Quan and Chang \cite{QC89,QC89a}, the DQMs have been employed with various types of base functions, among others, cubic B-spline DQM \cite{KD13,KD13a}, modified cubic B-spline differential quadrature method (MCB-DQM) \cite{AS13,MB14},  DQM based on fourier expansion and Harmonic  function \cite{SC97,SX97,SWB95}, sinc DQM \cite{KD11a}, generalized DQM \cite{SR92}, polynomial based DQM \cite{K10,JPM12}, quartic B-spline based DQM \cite{BKG15}, Quartic and quintic B-spline methods \cite{KD16}, exponential cubic B-spline based DQM \cite{KA15}, extended cubic B-spline DQM \cite{KA16}.

In this paper our aim is to develop a new method so called modified trigonometric cubic-B-spline differential quadrature method (MTB-DQM) for hyperbolic partial differential equations. Specially, MTB-DQM is employed for the numerical computation of two dimensional second order linear hyperbolic telegraph equation with both Dirichlet boundary conditions and Neumann boundary conditions. The MTB-DQM is a new differential quadrature method based on modified trigonometric cubic-B-splines as base functions. The MTB-DQM converts the initial-boundary value system of the telegraph equation into an initial value system of ODEs, in time. The resulting system of ODEs can be solved by using various time integration algorithm, among them, we prefer SSP-RK54 scheme (and SSP-RK43 scheme) \cite{GKS09,SR02} due to their reduce storage space, which results in less accumulation errors. The accuracy and adaptability of the method is illustrated by six test problems of second order linear hyperbolic telegraphic equations in dimension $(2+1)$.

The rest of the paper is organized into five more sections, which follow this introduction. Specifically, Section \ref{sec-disc} deals with
the description of MTB-DQM. Section \ref{sec-imp} is devoted to the procedure for the implementation of MTB-DQM for the initial value problem \eqref{int-tel-eqn1} with the boundary conditions \eqref{eqn-DBC} and \eqref{eqn-NBC}. The stability analysis of the MTB-DQM is studied in Section \ref{sec-stab}. Section \ref{sec-numer} is concerned with the main aim, the numerical study of six test problems, to establish the accuracy of the proposed method in terms of the relative error norm ($R_e$), $L_2$ and $L_\infty$ error norms, considering various grid points. Finally, Section \ref{sec-conclu} concludes the paper with reference to critical analysis and research perspectives.

\section{Description of MTB-DQM} \label{sec-disc}
The differential quadrature method is an approximation to derivatives of a function is the weighted sum of the functional values at certain nodes \cite{BKC72}. The weighting coefficients of the derivatives is depend only on grids \cite{SR92}. This is the reason for taking the partitions $P[\Omega]$ of the problem domain $\Omega=\{ (x, y)\in R^2: 0\leq x, y\leq 1\}$ distributed uniformly as follows:
 $$P[\Omega] = \{ (x_i, y_j)\in \Omega: h_x=x_{i+1}-x_{i}, h_y=y_{j+1}-y_{j}, i\in \Delta_x, j\in \Delta_y\},$$
where $\Delta_x=\{1,2,\ldots,N_x\}, \Delta_y=\{1,2,\ldots,N_y\}$, and $
h_x=\frac{1}{N_x-1} \mbox{ and } h_y=\frac{1}{N_y-1}$
are the discretization steps in both $x$ and $y$ directions, respectively. That is, a uniform partition in each $x, y$-direction with the grid points: $0=x_1 < x_2< \ldots <x_i<\ldots <x_{N_x-1}< x_{N_x}=1; ~0 = y_1 < y_2< \ldots< y_j< \ldots< y_{N_y-  1}<y_{N_y}=1.$

For generic grid point $(x_i, y_j)$, set $u_{ij}\equiv u_{ij}(t)\equiv u(x_i, y_j, t),~~i\in \Delta_x, j\in \Delta_y.$ Then $r$-th order derivative of $u(x,y,t)$, for $r \in \{1, 2\}$, with respect to $x, y$ at $(x_i, y_j)$ for $i\in \Delta_{x}, j\in \Delta_{y}$ is approximated as follows:
\begin{equation}\label{eq-deri1}
\begin{split}
& \left.\frac{\partial^r u}{\partial x^r}\right)_{ij} = \sum_{\ell=1}^{N_x} a_{i\ell}^{(r)} u_{\ell j},~~~  i \in \Delta_x; \qquad  \left.\frac{\partial^r u}{\partial y^r}\right)_{ij} = \sum_{\ell=1}^{N_y} b_{j \ell}^{(r)} u_{i\ell}, ~~~  j\in\Delta_y,
\end{split}
\end{equation}
%\end{flushleft}
where the time dependent unknown quantities $a_{i\ell}^{(r)}, ~~b_{j\ell}^{(r)}$ are the weighting functions of the $r$th-order derivatives, to be determine by adopting various type of base functions.

Trigonometric cubic B-spline function $T_i=T_i(x)$ at node $i$ in $x$ direction \cite{GV16,AMI14} read as:
\begin{eqnarray}\label{eq-cbs}
T_i= \frac{1}{\omega} \left\{ \begin{array}{ll}
                 p^3(x_i),  & x \in [x_{i}, x_{i+1}) \\
                 p(x_i)\{p(x_i)q(x_{i+2})+p(x_{i+1})q(x_{i+3})\} + p^2(x_{i+1})q(x_{i+4}), &   x \in [x_{i+1}, x_{i+2})\\
                 q(x_{i+4})\{p(x_{i+1}q(x_{i+3})+ p(x_{i+2})q(x_{i+4})\} + p(x_{i})q^2(x_{i+3}),
                    &  x \in [x_{i+2}, x_{i + 3})\\
                  q^3(x_{i+4}),  &  x \in [x_{i + 3}, x_{i + 4})
       \end{array}  \right.
\end{eqnarray}
where $p(x_i)=\sin(\frac{x-x_i}{2});~~ q(x_i)=\sin(\frac{x_i-x}{2})$, $\omega = \sin\(\frac{h_x}{2}\) \sin\(h_x\) \sin\(\frac{3h_x}{2}\)$. If we take,  $T_i=\zeta_{i+2}$, then the set $\{\zeta_{0}, \ldots, \zeta_{N_x}, \zeta_{N_x+1}\}$ forms a basis over the interval $[a, b]$. Setting
 \begin{equation*}
 \begin{split}
 & a_1 = \frac{\sin^2\(\frac{h_x}{2}\)}{ \sin(h_x) \sin\(\frac{3 h_x}{2}\)}; \quad a_2 = \frac{2}{ 1+ 2\cos(h_x)}; \quad a_3 = \frac{3}{ 4 \sin\(\frac{3 h_x}{2}\)};\\
 & a_4 = \frac{3 + 9 \cos(h_x)}{ 16 \sin^2\(\frac{h_x}{2}\) \(2 \cos\(\frac{h_x}{2}\) +  \cos\(\frac{3 h_x}{2}\)\) }; \quad a_5 =  \frac{ 3\cos^2(\frac{h_x}{2})}{ \sin^2\(\frac{h_x}{2}\) \(2 + 4 \cos\(h_x\)\) }.
 \end{split}
\end{equation*}

The values of $\zeta_i$ and its first and second derivatives in the grid point $x_j$, denoted by $\zeta_{ij}:=\zeta_i(x_j)$, $\zeta'_{ij}:=\zeta'_i(x_j)$ and $\zeta''_{ij}:=\zeta''_i(x_j)$, respectively, read:
\begin{eqnarray}\label{tab-coeff}
       \zeta_{ij}= \left\{ \begin{array}{ll}
                 a_2,   &\mbox{ if }  i-j=0 \\
                 a_1,  & \mbox{ if } i-j=\pm 1\\
                 0,  & \mbox {otherwise}
       \end{array}\right.; \quad
              \zeta'_{ij}= \left\{ \begin{array}{ll}
                  \pm a_3,  & \mbox{ if }  i-j=\pm1\\
                 0, &\mbox {otherwise}
       \end{array}\right.; \mbox{\quad}
       \zeta''_{ij}= \left\{ \begin{array}{ll}
                a_5,  &\mbox{ if }  i-j=0 \\
                a_4,  &\mbox{ if }  i-j=\pm 1\\
                 0  &  \mbox {otherwise}
       \end{array}\right.
\end{eqnarray}

The modified trigonometric cubic B-spline function are obtained by modifying the exponential cubic B-spline function \eqref{eq-cbs} as follows \cite{AS13}:
\begin{flushleft}
\begin{equation}\label{eq-modi-CBs}
\left\{ \begin{split}
& \psi_1 (x) = \zeta_1(x) + 2\zeta_0(x)\\
&  \psi_2 (x) = \zeta_2(x)  - \zeta_0(x)\\
& \vdots\\
& \psi_j (x) = \zeta_j(x), \mbox{ for } j = 3,4, \ldots, N_x-2\\
& \vdots\\
&  \psi_{N_x-1}(x) = \zeta_{N_x-1}(x)  - \zeta_{N_x+1}(x)\\
&  \psi_{N_x}(x) = \zeta_{N_x}(x)  + 2 \zeta_{N_x+1}(x)
\end{split}\right.
\end{equation}
\end{flushleft}

Now, the set $\{\psi_1, \psi_2,\ldots, \psi_{N_x}\}$ be a basis over $[a, b]$. The procedure for defining the  modified trigonometric cubic B-splines in $y$ direction is followed analogously.

\subsection{\textbf{The evaluation of the weighting coefficients $a_{ij}^{(r)}$ and  $ b_{ij}^{(r)} (r=1,2)$}}
In order to evaluate the weighting coefficients $a_{i\ell}^{(1)}$ of first order partial derivative in Eq. \eqref{eq-deri1}, the modified trigonometric cubic B-spline $\psi_p(x)$, $p\in\Delta_x$ are used in DQM as base functions. Setting $\psi'_{pi}:=\psi_p'(x_i)$ and $\psi_{p\ell}:=\psi_p(x_\ell)$. The approximate values of the first-order derivative is obtained using TBM-DQM as follows:
\begin{equation}\label{eq-deri1approxx}
\psi'_{pi} = \sum_{\ell\in \Delta_x} a_{i\ell}^{(1)}  \psi_{p\ell}, \qquad  p,i \in \Delta_x.
\end{equation}

Setting  $\Psi=[\psi_{p\ell}]$, $A=[a_{i\ell}^{(1)}]$ (the unknown weighting coefficient matrix), and $\Psi'=[\psi'_{p\ell}]$, then Eq. \eqref{eq-deri1approxx} can be re-written as the following set of system of linear equations:
\begin{equation}\label{eq-tri_syst}
\Psi A^T=\Psi',
\end{equation}
where the coefficient matrix $\Psi$ of order $N_x$ can be obtained from  \eqref{tab-coeff} and \eqref{eq-modi-CBs}:
\begin{equation*}
 \Psi= \left[
  \begin{array}{cccccccc}
  a_2+2a_1  $\quad$ & a_1 $\quad$&    $\quad$&           $\quad$&           $\quad$&         $\quad$&         \\
  $0    \quad$& a_2 $\quad$& a_1          $\quad$&           $\quad$&           $\quad$&         $\quad$&         \\
       $\quad$& a_1   $\quad$& a_2       $\quad$& a_1     $\quad$&           $\quad$&         $\quad$&         \\
       $\quad$&       $\quad$&  \ddots      $\quad$&   \ddots  $\quad$&   \ddots  $\quad$&         $\quad$&         \\
       $\quad$&       $\quad$&              $\quad$&   a_1     $\quad$&  a_2 $\quad$&    a_1  $\quad$&          \\
       $\quad$&       $\quad$&              $\quad$&           $\quad$&   a_1     $\quad$&   a_2$\quad$&    $0$   \\
       $\quad$&       $\quad$&              $\quad$&           $\quad$&           $\quad$&   a_1   $\quad$&   a_2+2a_1   \\
       \end{array}
\right]
\end{equation*}
and in particular the columns of the matrix $\Psi'$ read:
\begin{equation*}
\Psi'[1] = \left[ \begin{array}{c}
    2a_3\\
    -2a_3\\
    0 \\
  \vdots  \\
          \\
  0     \\
  0     \\
\end{array}
\right],
\Psi'[2] =\left[\begin{array}{c}
     a_3 \\
    $0$    \\
     -a_3  \\
     $0$   \\
   \vdots  \\
           \\
    $0$    \\
\end{array}\right],
\ldots,
\Psi'[N_x-1] =  \left[ \begin{array}{c}
     $0$     \\
     \vdots  \\
             \\
     $0$     \\
     a_3  \\
     $0$     \\
     -a_3  \\
\end{array}
\right], \mbox{ and }
\Psi'[N_x] =\left[ \begin{array}{c}
    $0$ \\
        \\
\vdots  \\
        \\
$0$     \\
  { }2a_3 \\
 -2a_3  \\
\end{array}
\right].
\end{equation*}
Notice that coefficient matrix $\Psi$ in Equation \eqref{eq-tri_syst} is diagonally dominant. For any $i\in \Delta_x$, the coefficients $a_{i1}^{(1)}, a_{i2}^{(1)}, \ldots, a_{iN_x}^{(1)}$ have been determined by adopting ``Thomas Algorithm" to solve system \eqref{eq-tri_syst}. The weighting coefficients $b_{i \ell}^{(1)}$ can be determined in similar manner, considering the grids in $y$ direction.

The weighting coefficients $a_{i\ell}^{(r)}$ and $b_{i\ell}^{(r)}$, for $r\ge 2$, can also be computed using the weighting functions in quadrature formula for second derivative on the given basis. But, in the present paper, we prefer the following recursive formulae \cite{SR92}:
\begin{flushleft}
 \begin{equation}\label{eq-coeff2}
\left\{ \begin{split}
& a_{i j}^{(r)} = r \left( a_{i j}^{(1)} a_{i i}^{(r-1)} - \frac{a_{i j}^{(r-1)}}{ x_i - x_j}\right), i\ne j: i, j\in \Delta_x,
\\&  a_{i i}^{(r)} = - \sum_{i = 1, i \ne j}^{N_x} a_{i j}^{(r)}, i= j: i, j\in \Delta_x.\\
& b_{i j}^{(r)} = r \left(b_{i j}^{(1)} b_{i i}^{(r-1)} - \frac{b_{i j}^{(r-1)}}{ y_i - y_j}\right), i\ne j: i, j \in \Delta_y
 \\& b_{i i}^{(r)} = - \sum_{i = 1, i \ne j}^{N_y} b_{i j}^{(r)}, i = j: i, j \in \Delta_y.\\
\end{split}\right.
\end{equation}
\end{flushleft}

\section{Implementation of MTB-DQM for the telegraph equation} \label{sec-imp}
 Keeping all above in mind, the initial value system \eqref{int-tel-eqn1} of second order telegraph equation under the transformation $u_t=v$ (and so, $u_{tt}=v_t$) transformed in to the initial value coupled system of first order differential equations as follows:
\begin{equation}\label{trans-tel-eqn2}
\left\{ \begin{array}{ll}
\frac{\partial u (x, y, t)}{\partial t}=v(x,y,t)\\ \\
\frac{\partial v (x, y, t)}{\partial t}= - 2\alpha  v (x, y, t) - \beta^2 u(x,y,t)+ \frac{\partial^2 u (x, y, t)}{\partial x^2} + \frac{\partial^2 u (x, y, t)}{\partial y^2}+f(x,y,t), (x,y)\in \Omega, t>0, \\
 u(x,y,0)=\phi(x,y), \qquad  v(x,y,0)=\psi(x,y), \qquad (x, y) \in \Omega.
\end{array}  \right.\end{equation}
Setting $f(x_i, y_j,t)=f_{ij}$. The MTB-DQM transforms equation \eqref{trans-tel-eqn2} to
\begin{equation}\label{trans-tel-eqn3}
\left\{ \begin{split}
& \frac{\partial u_{ij}}{\partial t}=v_{ij}\\
& \frac{\partial v_{ij}}{\partial t}= \sum_{\ell\in\Delta_x}a_{i\ell}^{(2)}u_{\ell j}+\sum_{\ell\in\Delta_y}b_{j\ell}^{(2)}u_{i\ell}-2\alpha v_{ij}-\beta^2 u_{ij}+f_{ij},\\
&  u_{ij}(t=0)=\phi_{ij}, \quad
 v_{ij}(t=0)=\psi_{ij}, \quad i \in \Delta_x, ~~  j\in \Delta_y.
\end{split}\right.
\end{equation}

Next, further simplification is not required in case of Dirichlet boundary conditions. In this case the solution on the boundary can be read directly from the conditions \eqref{eqn-DBC} as:
\begin{flushleft}
 \begin{equation}\label{eq-DBC-dis}
\left. \begin{split}
& u_{1j} = \phi_{1}(y_j, t)=\phi_{1}(j), \qquad u_{N_xj} = \phi_{2}(y_j, t)=\phi_{2}(j), \qquad  j\in \Delta_y, \\
& u_{i1} = \phi_{3}(x_i, t)=\phi_{3}(i), \qquad u_{iN_y} = \phi_{4}(x_i, t)=\phi_{4}(i), \qquad i \in \Delta_x,
\end{split}  \right\} t\geq0.
\end{equation}
\end{flushleft}
But, for Neumann or mixed type boundary conditions, further simplification is required. In this case, the solutions at the boundary are obtained by using MTB-DQM on the boundary, which yields a system of linear equations. On solving resulting system of the linear equation, the desired solution at the boundaries, is obtained. The procedure for implementation of Neumann boundary conditions can be seen in \cite{MB14}, and is given below to complete the problem.

Now, the functions $\psi_1, \psi_2$ from the boundary conditions \eqref{eqn-NBC} for $r=1$ yields the linear system in $u_{1j}, u_{N_x j}$, in compact matrix form as
\begin{equation}\label{eq-N2-dis}
\left[\begin{array}{cc}
a_{11}^{(1)}  & a_{1N_x}^{(1)}  \\
a_{N_x1}^{(1)} & a_{N_x N_x}^{(1)} \\
\end{array} \right] \left[
\begin{array}{c}
u_{1j} \\
u_{N_xj} \\
\end{array} \right]
=
\left[
\begin{array}{c}
S_1(j)\\
S_2(j) \\
\end{array} \right],
\end{equation}
Equation \eqref{eq-N2-dis} yields the solutions at the boundary of $x$-axis
\begin{equation}\label{eq-N3-dis}
u_{1j} = \frac{S_1(j) a_{N_xN_x}^{(1)} - S_2(j) a_{1N_x}^{(1)} }{a_{11}^{(1)}a_{N_xN_x}^{(1)}-a_{N_x 1}^{(1)}a_{1N_x}^{(1)}},  \qquad
u_{N_x j} = \frac{S_2(j) a_{11}^{(1)} - S_1(j) a_{N_x 1}^{(1)}}{a_{11}^{(1)}a_{N_xN_x}^{(1)}-a_{N_x 1}^{(1)}a_{1N_x}^{(1)}}, \qquad j\in \Delta_y.
\end{equation}
where $S_1(j)= \psi_1(j)-\sum_{\ell =2}^{N_x-1} a_{1\ell}^{(1)} u_{\ell j}$ and $ S_2(j)= \psi_1(j)- \sum_{\ell=2}^{N_x-1} a_{N_x \ell}^{(1)} u_{\ell j}$.

Similarly, the linear system obtained for $\psi_3, \psi_4$ from the boundary condition \eqref{eqn-NBC} yields the following
\begin{equation}\label{eq-N2y-dis}
%\left\{
 \begin{split}
&u_{i1} = \frac{S_3(i) b_{N_yN_y}^{(1)} - S_4(i) b_{1N_y}^{(1)} }{b_{11}^{(1)}b_{N_y N_y}^{(1)}-b_{N_y 1}^{(1)}b_{1N_y}^{(1)} }, \qquad  u_{iN_y} = \frac{S_4(i) b_{11}^{(1)} - S_3(i) b_{N_y 1}^{(1)}}{b_{11}^{(1)}b_{N_y N_y}^{(1)}-b_{N_y 1}^{(1)}b_{1 N_y}^{(1)} },
\end{split}%\right.
\qquad i\in \Delta_x,
\end{equation}
where $S_3(i)= \psi_3(i) - \sum_{\ell=2}^{N_y-1} b_{1\ell }^{(1)} u_{i\ell }$ and $ S_4(i)= \psi_4(i) - \sum_{\ell=2}^{N_y-1} b_{N_y \ell}^{(1)} u_{i\ell}$.

Finally, on using boundary values $u_{1j}, u_{N_x j}, u_{i1}$ and $u_{iN_y}$ obtained from either boundary conditions (Dirichlet boundary conditions \eqref{eqn-DBC} or Neumann boundary conditions  \eqref{eqn-NBC}), Eq. \eqref{trans-tel-eqn3} can be rewritten as follows:
\begin{equation}\label{eq-ode-finl1}
 \left\{ \begin{split}
& \frac{\partial u_{ij}}{\partial t}=v_{ij}\\
& \frac{\partial v_{ij}}{\partial t}= \sum_{\ell=2}^{N_x-1}a_{i\ell}^{(2)}u_{\ell j}+\sum_{\ell=2}^{N_y-1}b_{j\ell}^{(2)}u_{i\ell}-2\alpha v_{ij}-\beta^2 u_{ij}+K_{ij},\\
&  u_{ij}(t=0)=\phi_{ij}, \quad
 v_{ij}(t=0)=\psi_{ij},
\end{split}\right.
\end{equation}
where $2\leq i\leq N_x-1, 2\leq j\leq N_y-1$ and
 \begin{equation}\label{eq-F}
 \begin{split}
&K_{ij}=f_{ij}+ a_{i1}^{(2)} u_{1j} + a_{i N_x}^{(2)} u_{N_x j} +b_{j1}^{(2)} u_{i1} + b_{j N_y}^{(2)}u_{iN_y} .
%H_{ij}&=K_{ij}+ a_{i1}^{(2)} \phi_{1}(j) + a_{i N_x}^{(2)} \phi_{2}(j) +b_{j1}^{(2)}\phi_{3}(i) + b_{j N_y}^{(2)}\phi_{4}(i).
\end{split}
\end{equation}
A lot time integration schemes have been proposed for the numerical computation of initial value system of differential equations \eqref{eq-ode-finl1}, among others, the SSP-RK scheme allows low storage and large region of absolute property \cite{GKS09,SR02}. In what follows, we adopt SSP-RK43 scheme (and SSP-RK54 scheme) to integrate system from $t_m$ to $t_m+dt$ as given below, for getting the stable solution of hyperbolic differential equations:

\noindent SSP-RK43 Scheme:
\begin{equation*}
\begin{split}
&k^{(1)}= u^m + \frac{\triangle t}{2} L(u^m); ~~k^{(2)}= k^{(1)} + \frac{\triangle t}{2} L(k^{(1)}); ~~k^{(3)}= \frac{2}{3} u^m + \frac{k^{(2)}}{3} + \frac{\triangle t}{6} L(k^{(2)})\\
&u^{m + 1}= k^{(3)} + \frac{\triangle t}{2} L(k^{(3)}),\\
\end{split}
\end{equation*}

\noindent SSP-RK54 Scheme:
\begin{equation*}
\begin{split}
&k^{(1)} = u^m+ 0.391752226571890 \triangle t L(u^m) \\
&k^{(2)} = 0.444370493651235  u^m+ 0.555629506348765 k^{(1)} + 0.368410593050371 \triangle t L(k^{(1)}) \\
&k^{(3)}= 0.620101851488403 u^m+ 0.379898148511597 k^{(2)} + 0.251891774271694 \triangle t L(k^{(2)}) \\
&k^{(4)}= 0.178079954393132 u^m+ 0.821920045606868 k^{(3)} + 0.544974750228521 \triangle t L(k^{(3)}) \\
&u^{m + 1}= 0.517231671970585 k^{(2)}+ 0.096059710526147 k^{(3)} \\
&~~ + 0.063692468666290 \triangle t L(k^{(3)}) + 0.386708617503269 k^{(4)}+ 0.226007483236906 \triangle t L(k^{(4)})
\end{split}
\end{equation*}
$L()$ be linear differential operator as defined in \eqref{eq-TEL-ode1}.

\section{Stability analysis} \label{sec-stab}
The system \eqref{eq-ode-finl1}, in compact form, can be rewritten as:
 \begin{equation}\label{eq-TEL-ode1}
 \left\{ \begin{split}
&\frac{dU}{dt} = AU  + G, \\
& U (t=0)=U_0
\end{split}\right.
\end{equation}
where
\begin{enumerate}
   \item [$1)$] $A =\left[\begin{array}{ccc}
     O \qquad &  &\qquad I \\
     B \qquad & &\qquad-2\alpha I\\
       \end{array} \right],$ \qquad $G = \left[\begin{array}{c}
     O_1 \\
     K \\
       \end{array} \right]$, $U = \left[\begin{array}{c}
     u \\
     v\\
       \end{array} \right], $ and $U_0 = \left[\begin{array}{c}
     \phi \\
     \psi\\
       \end{array} \right]$
   \item [$1)$] $O$ and $O_1$ are null matrices;
   \item [$2)$] $I$ is the identity matrix of order $(N_x-2)(N_y-2)$;
   \item [$3)$] $U=(u,v)^T$ the vector solution at the grid points:

 \noindent  $u=(u_{22},u_{23}, \ldots,u_{2(N_y-1)}, u_{32},u_{33}, \ldots,u_{3(N_y-1)}, \ldots, u_{(N_x-1)2}, \ldots,u_{(N_x-1)(N_y-1)})$.

 \noindent  $v=(v_{22},v_{23}, \ldots,v_{2(N_y-1)}, v_{32},v_{33}, \ldots,v_{3(N_y-1)},  \ldots, v_{(N_x-1)2}, \ldots,v_{(N_x-1)(N_y-1)})$.

   \item [$4)$] $K=(K_{22},K_{23}, \ldots,K_{2(N_y-1)}, K_{32}, \ldots,K_{3(N_y-1)}, \ldots K_{(N_x-1) 2}, \ldots K_{(N_x-1) (N_y-1)} $, where $K_{ij}$, for $i \in \Delta_x, j \in \Delta_y$  is calculated from Eq. \eqref{eq-F}.
  \item [$5)$] $B = -\beta^2 I + B_x + B_y$, where $B_{x}$ and $B_{y}$ are the following matrices (of order $(N_x-2)(N_y-2)$) of the weighting coefficients $a_{ij}^{(2)}$ and $b_{ij}^{(2)}$:
\begin{equation}
\begin{array}{ll}
 B_x = \left[
    \begin{array}{cccc}
      a_{22}^{(2)} I_x & a_{23}^{(2)} I_x & \ldots & a_{2(N_x-1)}^{(2)} I_x \\
      a_{32}^{(2)} I_x & a_{33}^{(2)} I_x &\ldots  & a_{3(N_x-1)}^{(2)} I_x \\
      \vdots & \vdots & \ddots & \vdots \\
      a_{(N_x-1)2}^{(2)} I_x & a_{(N_x-2)3}^{(2)} I_x & \ldots  & a_{(N_x-1)(N_x-1)}^{(2)} I_x \\
    \end{array}
  \right], ~~ & \begin{array}{ll}
   B_y = \left[
    \begin{array}{cccc}
      M_y & O_y & \ldots & O_y \\
      O_y & M_y &\ldots  & O_y \\
      \vdots & \vdots & \ddots & \vdots \\
      O_y & O_y & \ldots  & M_y \\
    \end{array}
  \right]
\end{array}
\end{array}
\end{equation}
where  identity matrix, $I_x$, and null matrix, $O_y$, both are of order $(N_y-2)$ and
\begin{equation*}
\begin{array}{ll}
M_y = \left[
    \begin{array}{cccc}
      b_{22}^{(2)} & b_{23}^{(2)} & \ldots & b_{2(N_y-1)}^{(2)} \\
      b_{32}^{(2)} & b_{33}^{(2)} & \ldots & b_{3(N_y-1)}^{(2)} \\
      \vdots & \vdots & \ddots & \vdots \\
      b_{(N_y-1)2}^{(2)} & b_{(N_y-1)3}^{(2)} & \ldots & b_{(N_y-1)(N_y-1)}^{(2)}
    \end{array}
  \right]
\end{array}
\end{equation*}
 \end{enumerate}
The stability of MTB-DQM for the telegraph equation \eqref{int-tel-eqn1} depends on the stability of the system of ODEs  defined in \eqref{eq-TEL-ode1} and the proposed MTB-DQM for temporal discretization may not converge to the exact solution whenever the system of ODEs \eqref{eq-TEL-ode1} is unstable \cite{Jain83}. The stability of \eqref{eq-TEL-ode1} depends on the eigenvalues of the matrix $A$ \cite{Jain83}.  In fact, the stability region is the set $\mathcal{S}= \{z \in C: \mid R(z)\mid \leq 1, z = \lambda_A \triangle t \}$, where $R(.)$ is the stability function and $\lambda_A$ is the eigenvalue of the matrix $A$. Specially, the stability region for SSP-RK43 scheme \cite[pp.72]{Jain83} and SSP-RK54 scheme  \cite[Fig.5]{KBYD}, which shows that the sufficient condition for the stability of the system \eqref{eq-TEL-ode1} is that for some $\triangle t$, $\lambda_A \triangle t \in \mathcal{S}$ for each eigenvalue $\lambda_A$ the matrix $A$, and so, the real part of each eigenvalue is necessarily less than or equal to zero.

Let $\lambda_{A}$ be an eigenvalue of $A$ associated with the eigenvector $(X_1,X_2)^T$, where both $X_1,X_2$ are vector of order $(N_x-2)(N_y-2)$, then from Eq. \eqref{eq-TEL-ode1}, we get
\begin{equation}\label{eq-egen}
A \left[\begin{array}{c}
     X_1 \\
     X_2\\
       \end{array} \right]=\left[\begin{array}{ccc}
     O \qquad & &\qquad  I \\
     B \qquad & &  \qquad  -2 \alpha I\\
       \end{array} \right]
\left[\begin{array}{c}
     X_1 \\
     X_2\\
       \end{array} \right] = \lambda_A \left[\begin{array}{c}
     X_1 \\
     X_2\\
       \end{array} \right],
\end{equation}

The above equation yields the following two equations
\begin{equation}\label{eqn-eign-a}
IX_2 = \lambda_A X_1,
\end{equation}
and
\begin{equation}\label{eqn-eign-b}
BX_1 - 2\alpha X_2 = \lambda_A X_2.
\end{equation}
Eq. \eqref{eqn-eign-a} and Eq. \eqref{eqn-eign-b} yields
%\begin{equation}\label{eq-eigenB}
$ BX_1 = \lambda_A (\lambda_A + 2\alpha) X_1$
%\end{equation}
, that is the eigenvalue $\lambda_B$ of the matrix $B$ is given by $\lambda_B=\lambda_A (\lambda_A + 2 \alpha)$, where the matrix $B$ is given by
\begin{equation}\label{eq-matrixB}
B= -\beta^2 I + B_x + B_y,
\end{equation}

\begin{figure}
\includegraphics[height=7.0cm, width=16.95cm]{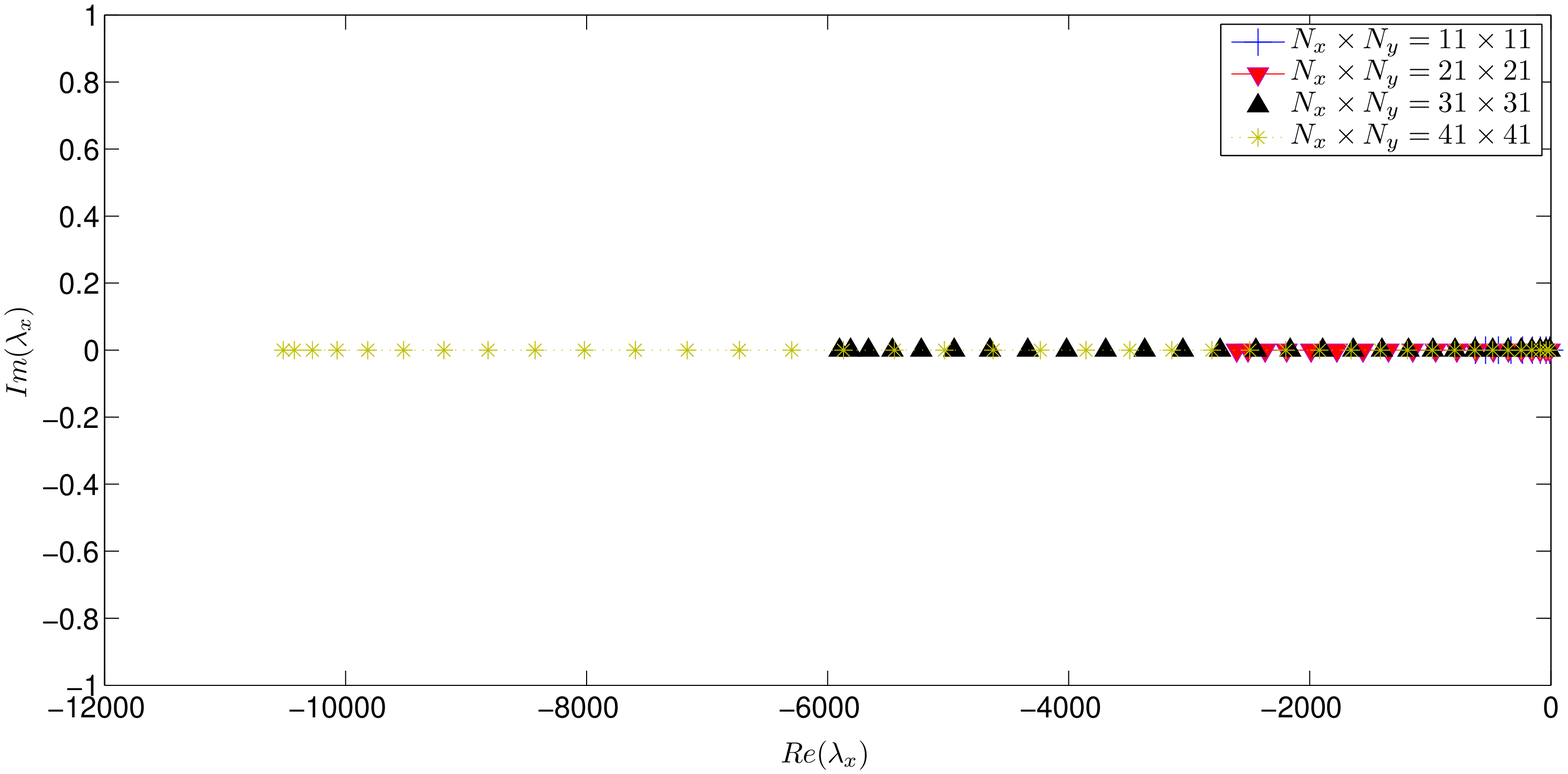}
\includegraphics[height=7.0cm, width=16.95cm]{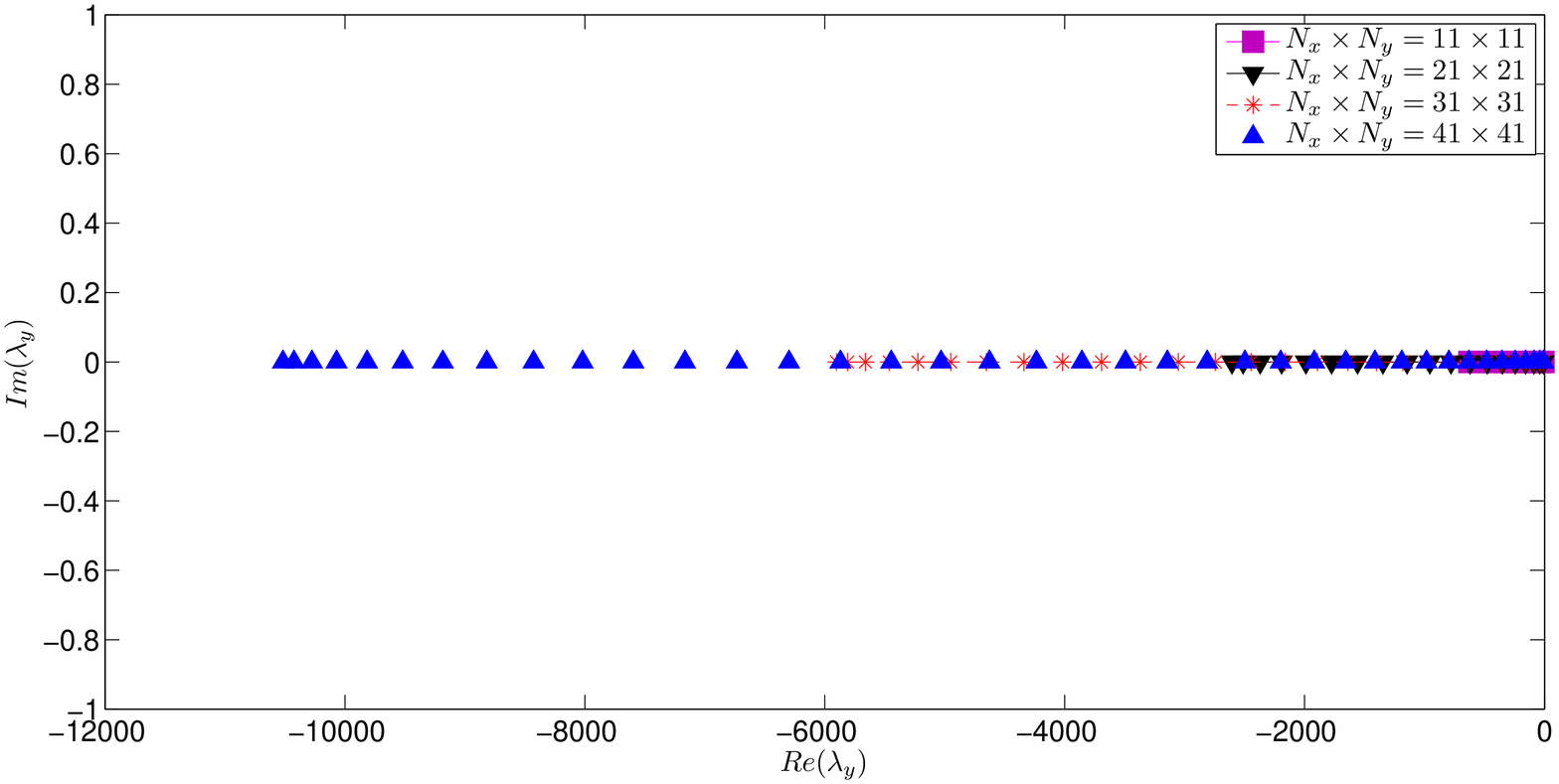}
\caption{ Eigenvalues $\lambda_x$ and $\lambda_y$ of $B_x$ and $B_y$ resp., for different grid sizes $N_x\times N_y = 11\times11, 21\times21, 31\times31, 41\times41$}
\label{FigB12}
\end{figure}

For different grid sizes: $N_x\times N_y = 11\times11, 21\times21, 31\times31, 41\times41$, the eigenvalues $\lambda_x$ and $\lambda_y$ have been depicted in Fig \ref{FigB12}. Eq. \eqref{eq-matrixB} and Fig. \ref{FigB12} confirms that for different values of the grid sizes the computed eigenvalues $\lambda_B=\lambda-\beta^2$ of $B$ are real negative numbers, that is
\begin{equation}\label{eq-A}
Re\left(\lambda_B\right) \leq 0 \mbox{~~  and ~~} Im\left(\lambda_B\right) = 0,
\end{equation}
where $Re(z)$ and $Im(z)$ denote the real and the imaginary part of $z$, respectively. Setting $\lambda_A = x+\iota y$, then
\begin{equation}\label{eq-B}
\begin{split}
\lambda_B &=\lambda_A (\lambda_A + 2\alpha)\\
& =x^2-y^2 + 2 \alpha x + 2 \iota (x+\alpha)y.
\end{split}
\end{equation}
Eq. \eqref{eq-A} and Eq. \eqref{eq-B} implies together that
\begin{equation}\label{eq-sta1}
 \left\{ \begin{split}
& x^2-y^2 + 2 \alpha x < 0\\
&(x+\alpha)y=0\\
\end{split}\right.
\end{equation}
Eq.\eqref{eq-sta1} have the following possible solutions
\begin{enumerate}
   \item  [$1)$] If $y\ne 0$, then $x=-\alpha,$
   \item  [$2)$] If $y=0$, then $(x+\alpha)^2 < \alpha^2.$
 \end{enumerate}
As $\alpha>0$, and so, in each case $x$ is negative. Thus, for a given $h$ one can find a sufficient small value of $\triangle t$ so that for each eigenvalue $\lambda_A$ of matrix $A$, $\triangle t \lambda_A$ belongs to the stability region $\mathcal{S}$, and so, the MTB-DQM produces stable solutions for the system of second order hyperbolic telegraph equations in $(2+1)$ dimension.

\section{Numerical experiments and discussion} \label{sec-numer}
This section deals with the main goal of the paper is the computation of numerical solution of the telegraph equation obtained by new algorithm: MTB-DQM, with space step size $h_x=h_y=h$. The accuracy and the efficiency of the method are measured for six test problems in terms of the discrete error norms: $L_2$, $L_{\infty}$ and relative error ($R_e$) \cite{DG10a} which are defined as:
\begin{equation*}
\begin{split}
& L_2 = \(h\sum_{k=1}^{N} e_{k}^2\)^{1/2}, \mbox{       } L_\infty = \max_{1\leq k\leq N} e_{k}, \mbox{   and   }
 R_e = \left(\frac{\sum_{k=1}^{N} e_{k}^2}{\sum_{k=1}^{N} u_{k}^2}\right)^{1/2}
\end{split}
\end{equation*}
where $e_{k}=\left|u_{k} - u_{k}^*\right|$; $u_{k}$ and $u_{k}^*$ denote exact solution and numeric solution at node $k$, respectively.

\begin{figure}
\centering
\includegraphics[height=6.50cm,width=6.26cm]{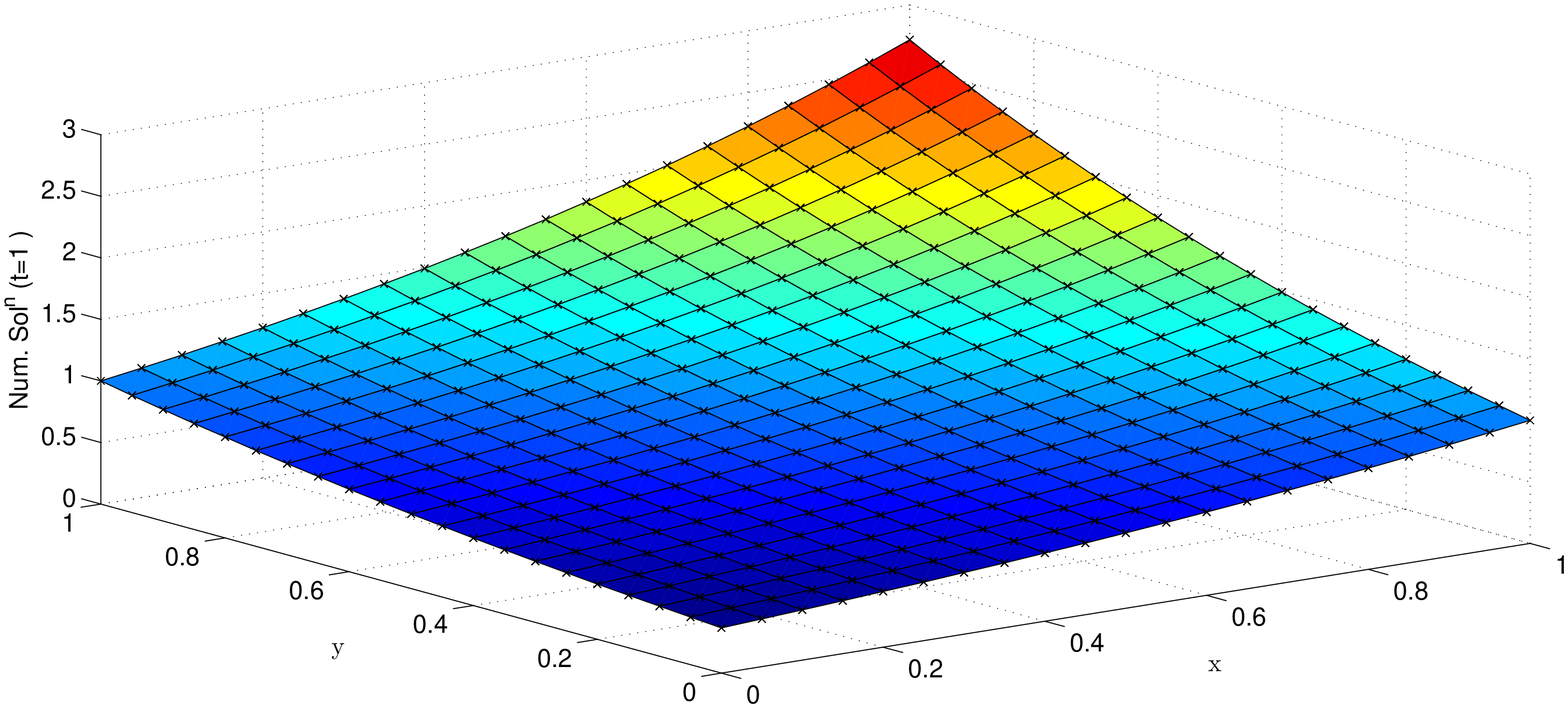}
\includegraphics[height=6.50cm,width=6.26cm]{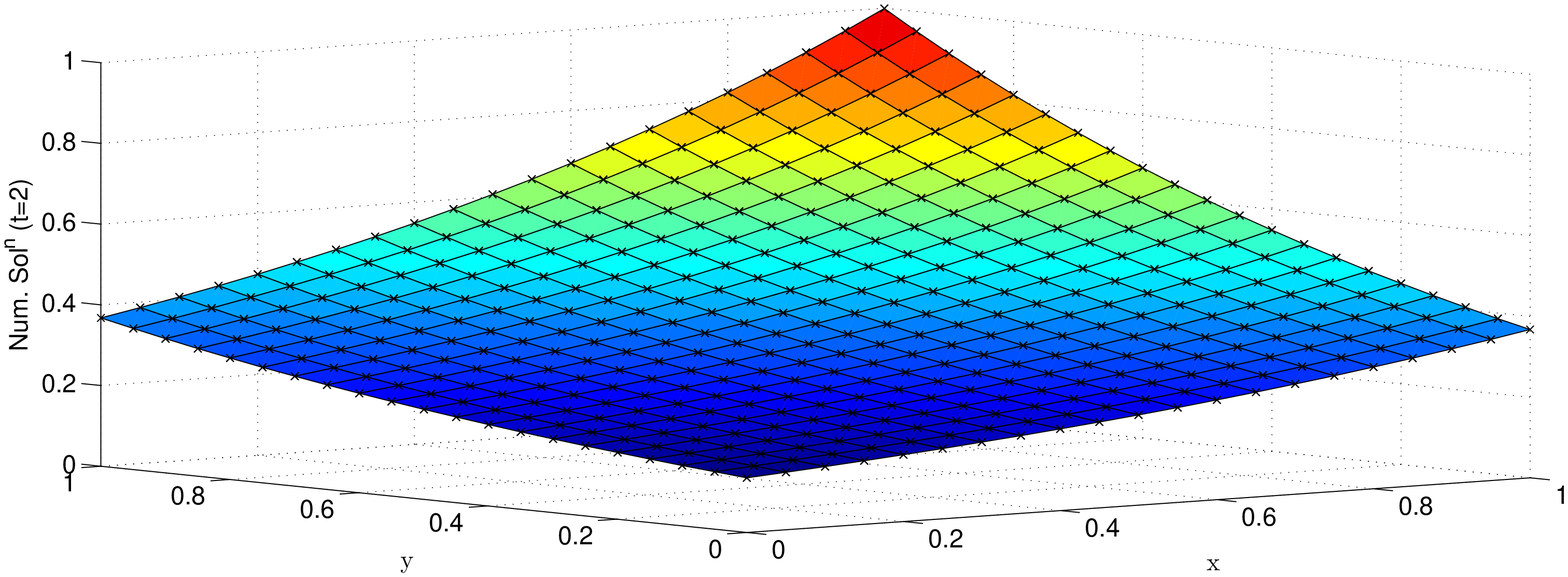}
\includegraphics[height=6.50cm,width=6.26cm]{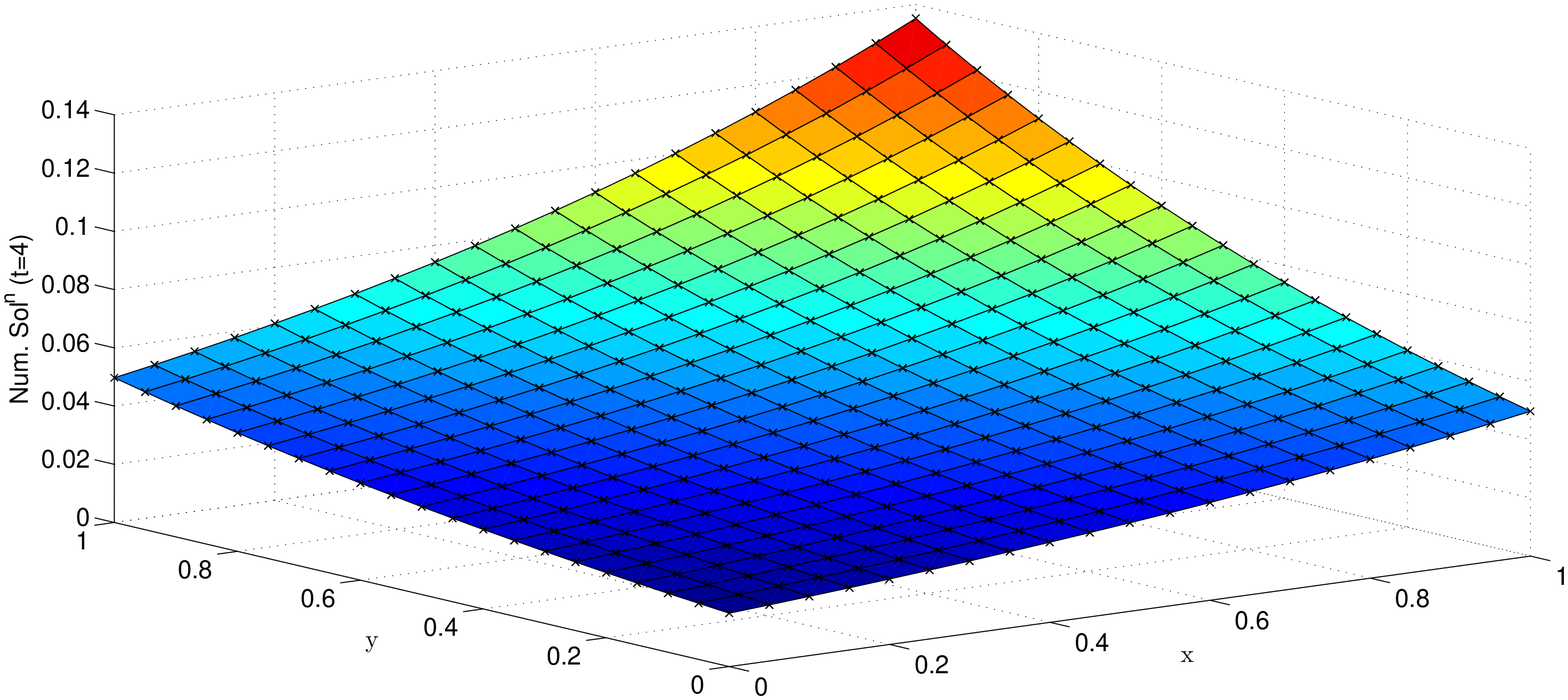}
\caption{Plots of MTB-DQM solutions at different time levels for Problem \ref{ex5}}\label{fig4.1}
\includegraphics[height=6.50cm,width=6.26cm]{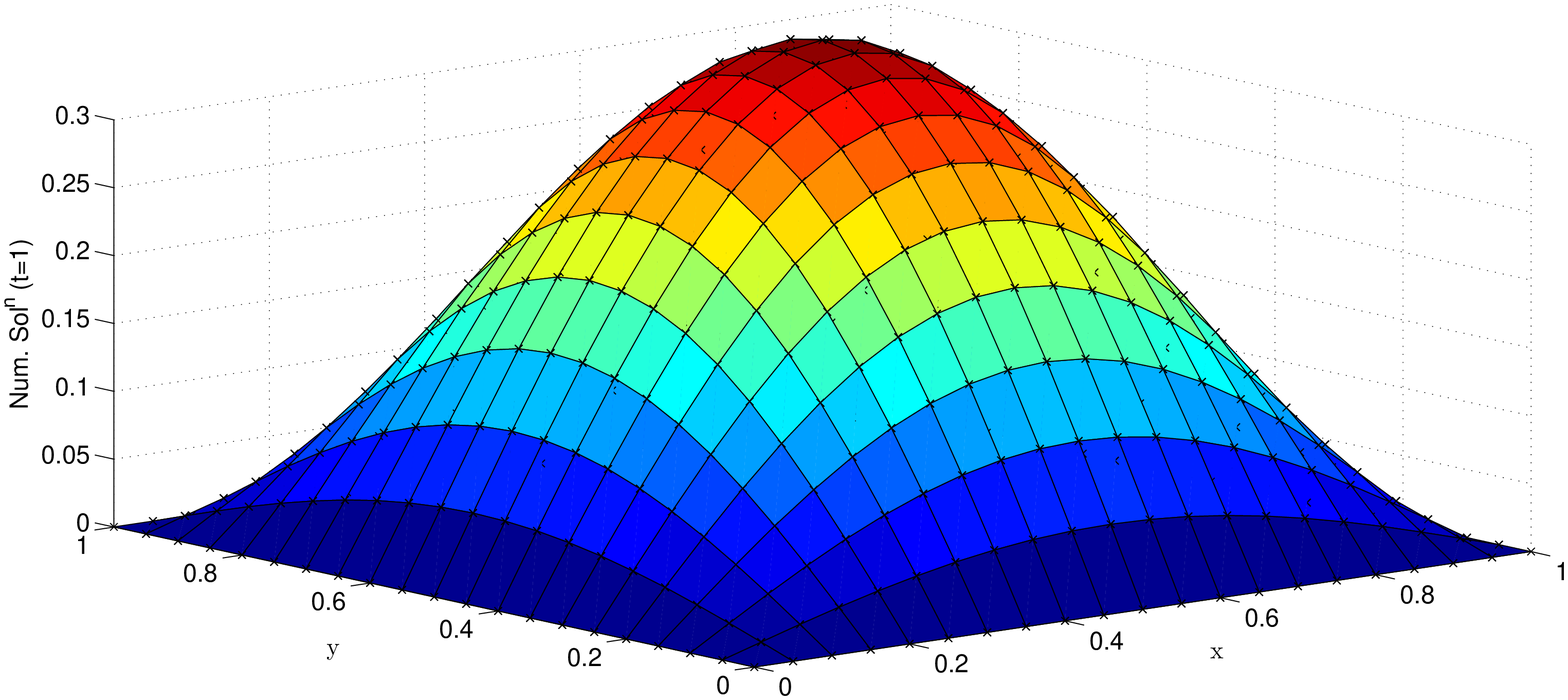}
\includegraphics[height=6.50cm,width=6.26cm]{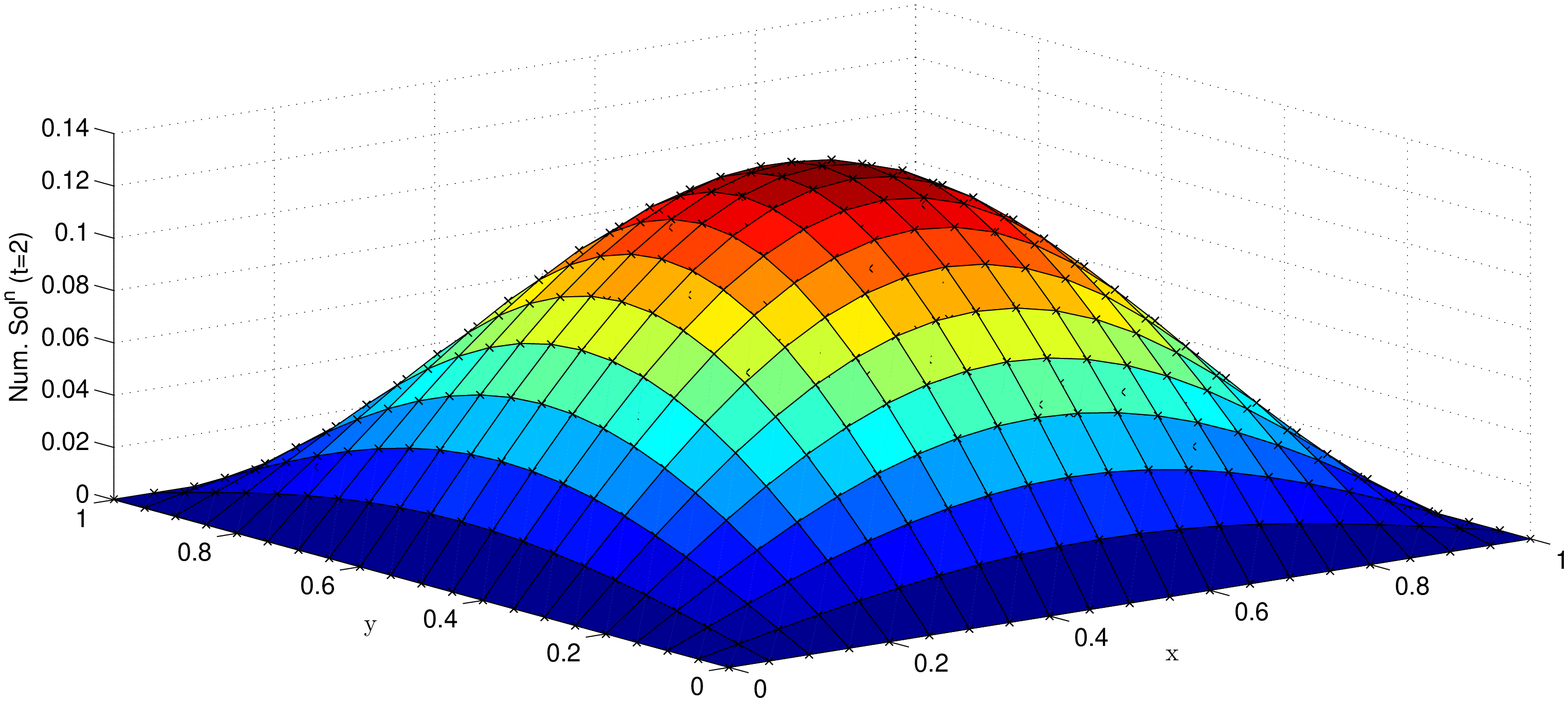}
\includegraphics[height=6.50cm,width=6.26cm]{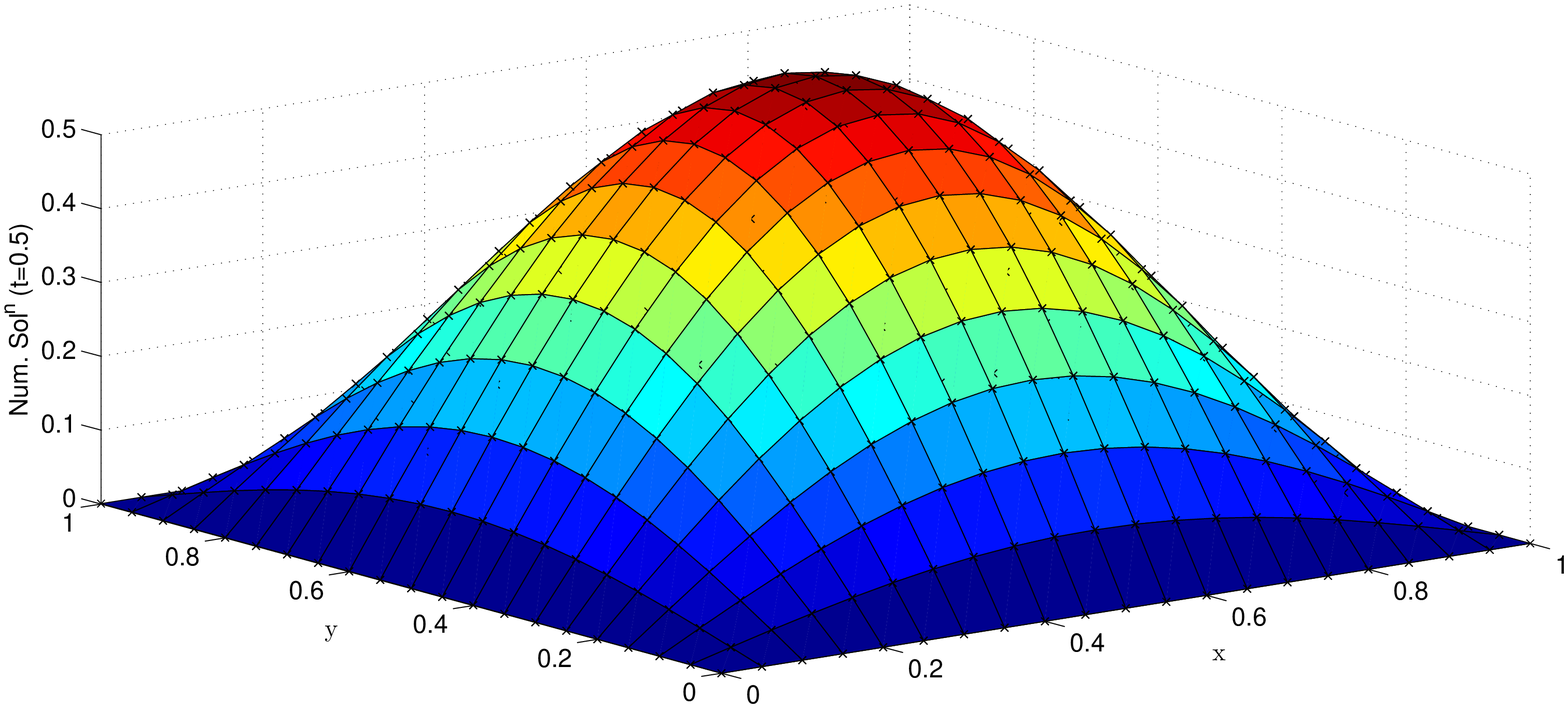}
\caption{Plots of MTB-DQM solutions at different time levels for Problem \ref{ex6}} \label{fig6.1}
\includegraphics[height=6.50cm,width=6.26cm]{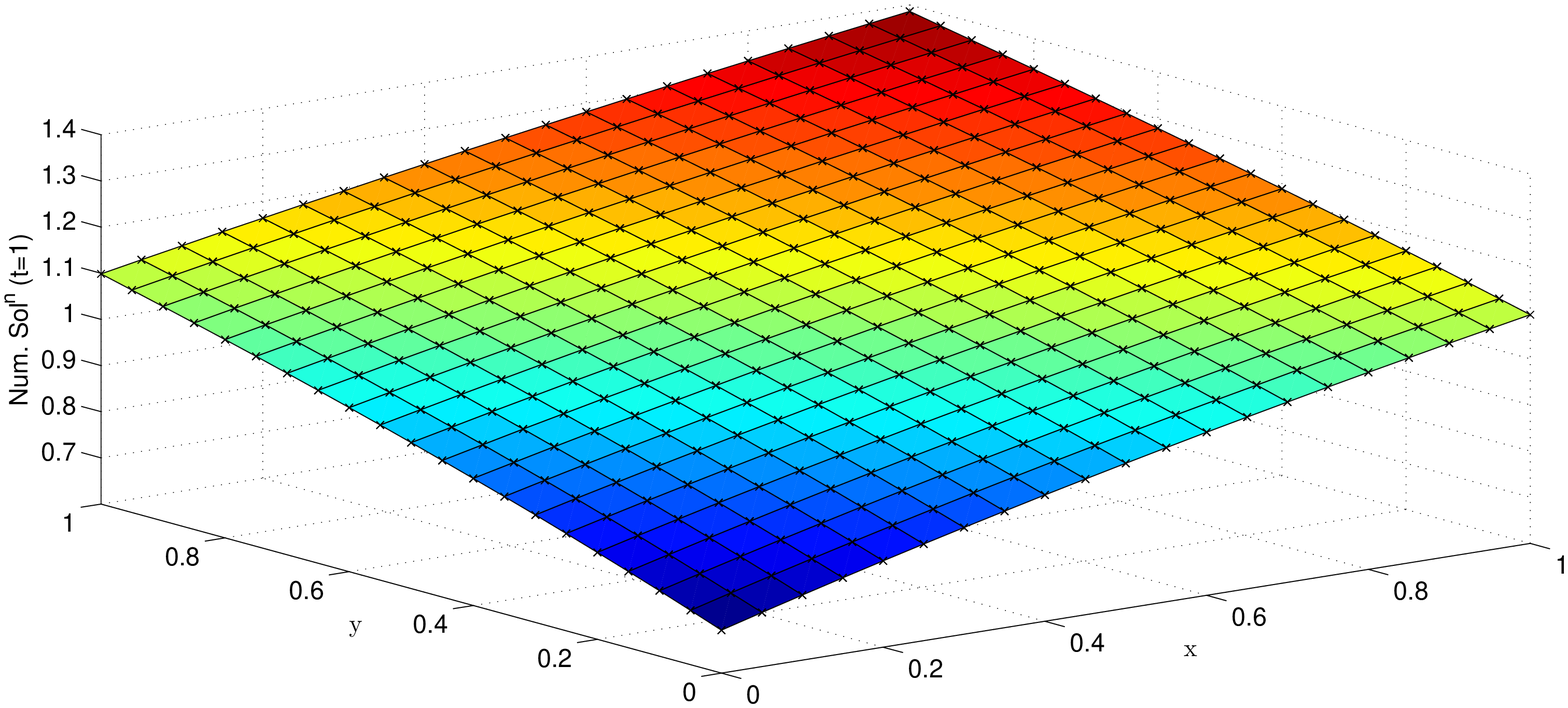}
\includegraphics[height=6.50cm,width=6.26cm]{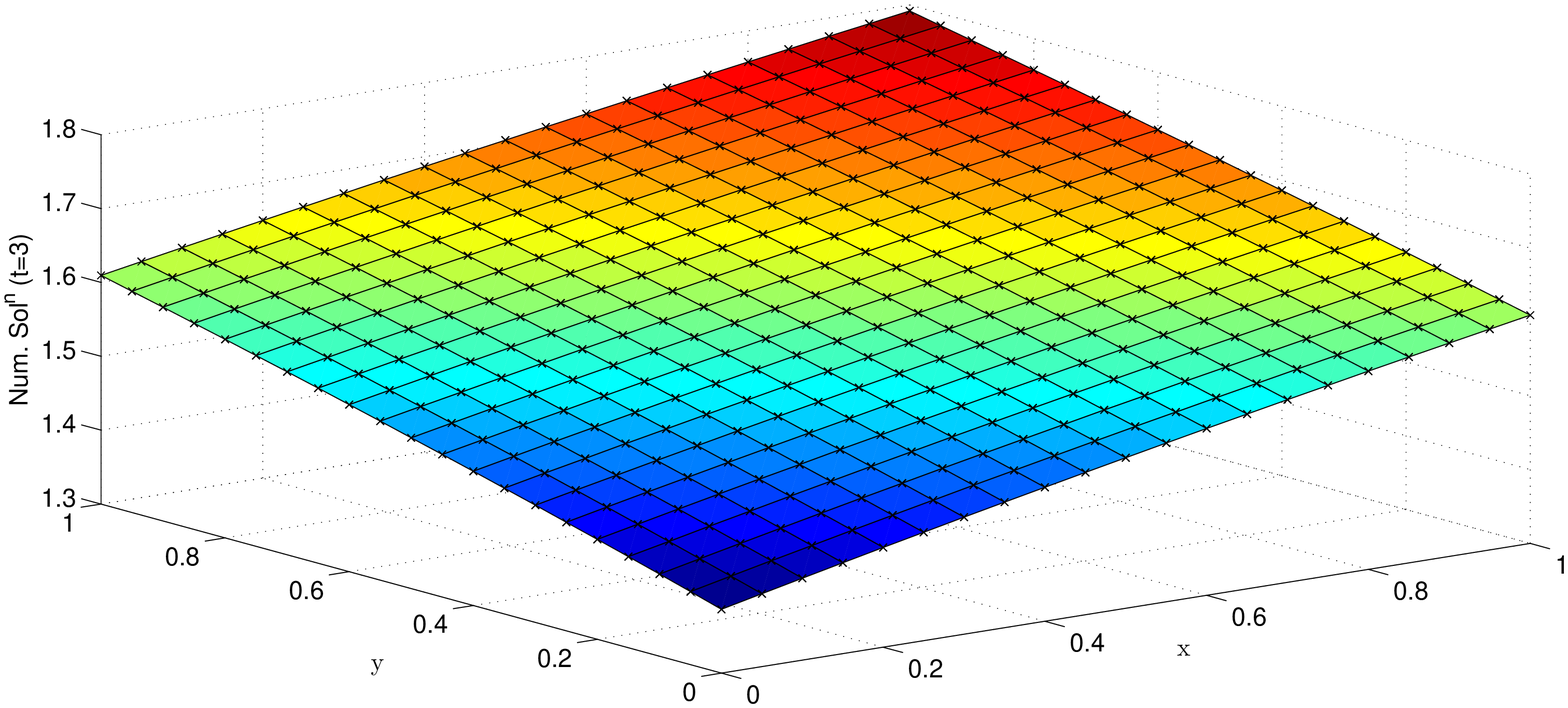}
\includegraphics[height=6.50cm,width=6.26cm]{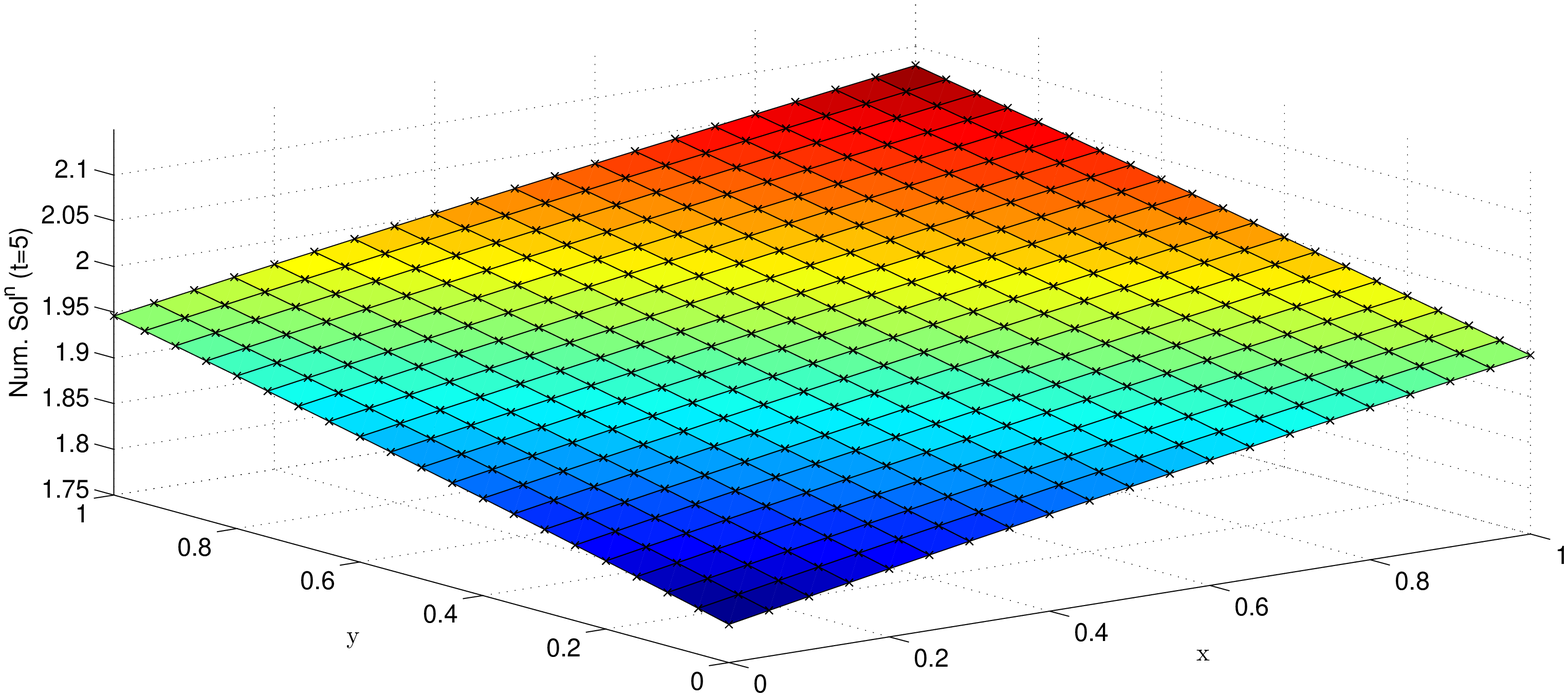}
\caption{Plots of MTB-DQM solutions at different time levels for Problem \ref{ex7}}\label{fig7.1}
\end{figure}

\begin{problem}\label{ex5}
Consider the telegraph equation \eqref{int-tel-eqn1} in the region $\Omega$ with $\alpha=1,\beta=1, f(x,y,t)=-2\exp{(x+y-t)}$,
and $ \phi(x,y)= \exp{(x+y)},\psi(x,y)=-\exp{(x+y)}$ in $\Omega$, and the mixed boundary conditions $
 \phi_1(y,t)=\exp {(y-t)}, \phi_2(y,t)=\exp {(1+y-t)} $ for $0\le y\le1 $ and $\psi_3(x, t)=\exp {(x-t)},
\phi_4(x,t)=\exp {(1+x-t)}$ for $ 0\le x\le 1.$
The exact solution \cite{DG10} is given by
\begin{equation}\label{ex5-exact-soln}
 u(x,y,t)=\exp{(x+y-t)}
\end{equation}
In Table \ref{tab5.1}, the error norms: $L_2, L_\infty, R_e$ have been compared with the errors obtained due to MCB-DQM \cite{MB14} taking the parameters $h=0.1, \triangle t=0.01$ and $h=0.05, \triangle t= 0.001$. The surface plots of the MTB-DQM solutions and the exact solutions at different time levels  $t=1, 2, 4$ have been depicted in Fig. \ref{fig4.1}. It is evident that the proposed solutions are more accurate in comparison to the results by Mittal and  Bhatia \cite{MB14}.
\end{problem}

\begin{problem}\label{ex6}
Consider the telegraph equation \eqref{int-tel-eqn1} with $\alpha=1, \beta=1$, $f(x,y,t)=2\pi^2 \exp(-t)\sin \pi x $ $\sin {\pi y}$  in region $\Omega, t>0 $ together with $ \phi(x,y)= \sin {\pi x}\sin {\pi y}, \psi(x,y)=-\sin {\pi x}\sin {\pi y}$ in $\Omega$, and the mixed boundary conditions
$ \psi_1(y,t)=\pi \exp {(-t)}\sin{(\pi y)}, \phi_2(y,t)=0$ in  $0 \le y\le 1,$ and
 $\phi_3(x,t)=0,  \psi_4(x,t)=-\pi \exp {(-t)}\sin{(\pi y)}$ in $ 0\le x\le 1$. The exact solution as given in \cite{DG10} is $ u(x,y,t)=\exp{(-t)}\sin{(\pi x)} \sin{(\pi y)}.$

The solutions are computed in terms of $L_2, L_\infty $ error norms, for $h=0.1, \triangle t =0.01$ and $h=0.05, \triangle t =0.001$, and are reported in Table \ref{tab-ex6.1}. The surface plots of MTB-DQM solutions and exact solutions at $t=0.5, 1, 2$ are depicted in Fig. \ref{fig6.1}. The findings from Table \ref{tab-ex6.1} concludes that the proposed MTB-DQM solutions are more accurate as compared to that obtained in \cite{MB14}.
\end{problem}

\begin{problem}\label{ex7}
Consider the telegraph equation \eqref{int-tel-eqn1} with the exact solution  $u(x,y,t)=\log(1+x+y+t)$ for $\alpha=\beta=1$ as given in \cite{DG10}, from which one can  extract the functions $\phi$ and $f$ directly. The mixed boundary conditions are taken as
$ \phi_1(y,t)=\log(1+y+t), \psi_2(y,t)=\frac{1}{2+y+t}$ for  $0\le y\le 1$ and $\psi_3(x,t)=\frac{1}{1+x+t}, \phi_4(x,t)=log(2+x+t)$ for $0\le x\le 1.$

The error norms: $L_2, L_\infty$ and $R_e$ has been computed with the parameters $\triangle t=0.001, h=0.05$ in the region $\Omega$ by adopting MTB-DQM (SSP-RK43), and are compared with the error norms due to MCB-DQM (SSP-RK43) \cite{MB14} and Dehghan and Ghesmati \cite{DG10a}, in Table \ref{tab-ex7.0}. It is seen from Table \ref{tab-ex7.0} that the computed MTB-DQM results are more accurate in comparison to the results obtained in \cite{MB14,DG10a}. The surface plots of the MTB-DQM solutions and the exact solutions at different time levels $t=1,2,3$ are depicted in Fig. \ref{fig7.1}.
\end{problem}

\begin{problem}\label{ex1}
Consider the telegraph equation \eqref{int-tel-eqn1} in the region $\Omega$ with $\alpha=\beta=1$, $f(x,y,t)=2(\cos t-\sin t)\sin x \sin y$,
$ \phi(x,y)=\sin x. \sin y; \psi(x,y)=0$, and the Dirichlet boundary conditions:
\begin{equation}\label{DC-ex1}
\left\{ \begin{array}{ll}
 \phi_1(y,t)=0, \phi_2(y,t)=\cos t \sin(1)\sin y, & 0\le y\le 1, \\ \\
  \phi_3(x,t)=0,   \phi_4(x,t)=\cos t \sin x \sin (1), \qquad& 0\le x \le 1,
\end{array}  \right.\end{equation}
The exact solution \cite{DG10} is
\begin{equation}\label{ex1-exact-soln}
 u(x,y,t)=\cos t \sin x \sin y
\end{equation}
We have compared the computed error norms: $R_e$, $L_2, L_\infty$  with the recent results of Mittal and Bhatia \cite{MB14} at different time levels $t\leq 10$, and reported them in Table \ref{tab-ex1.1} for the parameters $\triangle t=0.01, h=0.10$ while for $\triangle t=0.001, h=0.05$ in Table \ref{tab-ex1.2}. %The comparison of computed physical solution behavior with the the exact solution behavior at $t=1,2,3$ is depicted in Fig. \ref{fig1.1} for the parameters $\triangle t=0.001, h=0.05$.
The findings shows that the proposed solution are much better than that in \cite{MB14}, and are in excellent agreement with the exact solutions.
\end{problem}

\begin{problem}\label{ex2}
Consider the telegraph equation  \eqref{int-tel-eqn1} with $f(x,y,t)=(-2\alpha+\beta^2-1)\exp(-t)\sinh x \sinh y$,
 $\phi(x,y)=\sinh x \sinh y, \psi(x,y)=-\sinh x \sinh y,$ in $\Omega$;$\phi_1(y,t)=0; \phi_2(y,t)=\exp{(-t)}\sinh(1)\sinh y $ for $0\le y\le1$ and $
 \psi_3(x, t)=0; \psi_4(x,t)=\exp {(-t)}\sinh x \sinh (1)$ for $ 0\le x\le 1$.

The exact solution \cite{JPM12} is
\begin{equation}\label{ex2-exact-soln}
 u(x,y,t)=\exp{(-t)}\sinh x \sinh y
\end{equation}
The error norms are computed for $\alpha=10,\beta=5$ and $\alpha=10,\beta=0$ with parameters $\triangle  t =0.01, 0.001$, $h=0.1, 0.05$. The computed $L_2, L_\infty$  error norms for different time levels $t\le 10$ are compared with that obtained by \cite{MB14} in Table \ref{tab-ex2.1} taking $\triangle  t=0.01, h=0.1$. In Table \ref{tab-ex2.2}, the computed error norms for $\triangle t=0.001$ and $h=0.05$ are compared with that by Mittal and Bhatia \cite{MB14} and Jiwari et al. \cite{JPM12}. The findings from the above tables confirms that the MTB-DQM results are more accurate in compression to that obtained by MCB-DQM \cite{MB14} and PDQM \cite{JPM12}. %The surface plots of numerical and exact solutions at$ t=1,2,3$ with $\delta t=0.001$ and $h=0.05$ are depicted in Fig. \ref{fig2.1}.
\end{problem}

\begin{problem}\label{ex4}
 Consider the telegraph equation \eqref{int-tel-eqn1} in the region $\Omega$ with
$f(x,y,t)=(-3 \cos t+2\alpha \sin t +\beta^2 \cos t)\sinh x \sinh y$, and
$ \phi(x,y)= \sinh x \sinh y, \psi(x,y)=0$ in $\Omega$, and $ \phi_1(y,t)=o, \phi_2(y,t)=\cos t \sinh (1)\sinh y $ for $0\le y\le 1$, and $\phi_3(x,t)=0, \psi_4(x,t)=\cos t \sinh x\sinh (1)$ for $0\le x\le 1.$

The exact solution \cite{JPM12} is
\begin{equation}\label{ex4-exact-soln}
 u(x,y,t)=\cos t \sinh x \sinh y
\end{equation}

The error norms: $L_2, L_\infty, R_e$ at different time levels are computed for the parameters $\alpha=10, \beta=5 $ and $\alpha=50, \beta=5$ taking $\triangle t=0.001, h=0.05$, and are reported in Table \ref{tab-ex4.1}. The findings from Table \ref{tab-ex4.1} confirm that our results are comparably more accurate than that obtained by MCB-DQM \cite{MB14} and PDQM \cite{JPM12}. %A comparison of numerical solutions with exact solutions for $t=1, 2, 3$ are depicted in Fig. \ref{fig3.1}.
\end{problem}

\section{Conclusion} \label{sec-conclu}
In this paper, a new differential quadrature method based on modified trigonometric cubic B-splines as ne base functions has been developed, and so, call it modified trigonometric cubic B-spline differential quadrature method (MTB-DQM) for solving initial value system of linear hyperbolic partial differential equations in $(2+1)$ dimension. Specially, we implemented this method along with SSP-RK43 (and SSP-RK54) scheme for time integration for solving the second order hyperbolic telegraph equation in $(2+1)$ dimension subject to initial conditions, and each type of boundary conditions: Dirichlet, Neumann, mixed boundary conditions.

The matrix stability analysis method has been adopted to discus the stability of the MTB-DQM for the linear initial value system of second order hyperbolic telegraph equations. For different values of grid sizes: $N_x\times N_y = 11\times11, 21\times21, 31\times31, 41\times41$, we have computed the eigenvalues $\lambda_x, \lambda_y$ to compute the eigenvalues $\lambda_B$ of the matrix $B$ for arbitrary $\beta^2$, which are found are real and negative. Using this property of $\lambda_B$, we found that each eigenvalue $\lambda_A$ of the matrix $A$ has negative real part, which confirms that stability of the system of second order hyperbolic telegraph equations in $(2+1)$ dimension.

The findings from Section \ref{sec-numer}, shows that the MTB-DQM solutions are found comparatively more accurate than that obtained in \cite{MB14,JPM12,DG10a}.

\section*{Acknowledgement}
\noindent %The authors would like to thank the anonymous referees for their time, effort, and extensive comments which improve the quality of the presentation of the paper.
To be inserted ...

\newpage

\section{List of Tables and Figures}

\begin{table}[!htbp]
\caption{Comparison of the error norms for Problem \ref{ex1} for $h=0.1,\triangle t=0.01$ }\label{tab-ex1.1}
\vspace{.2cm}
\centering
\begin{tabular}{*9l*3l}
\toprule
{}& \multicolumn{3}{l}{MCB-DQM (SSP-RK43) \cite{MB14} }  &{}& \multicolumn{7}{l}{MTB-DQM}   \\ \cline{2-4} \cline{6-12}
{}&\multicolumn{3}{c}{  }  &{}& \multicolumn{3}{l}{SSP-RK54 } & & \multicolumn{3}{l}{SSP-RK43} \\ \cline{2-4} \cline{6-8} \cline{10-12}
{t} &$L_2$ &$L_\infty$& $R_e$ &{}&$L_2$ &$L_\infty$ &$R_e$ &&$L_2$ &$L_\infty$ &$R_e$\\
\midrule				
1	& 9.9722E-04\quad&2.2746E-03 \quad&5.9762E-03     &&3.7238E-06    \quad&4.5385E-06\quad&    2.4548E-05&	&    4.0186E-06\quad& 4.7532E-06  \quad& 2.6490E-05 \\
2	&1.0926E-03	&2.8706E-03	&8.5019E-03		&&4.4905E-06	&5.6389E-06&	3.8433E-05&	&	2.2235E-06&	3.6462E-06&	1.9030E-05	\\
3	&2.2877E-04	&6.0818E-04	&7.4720E-04		&&3.7914E-06	&6.3660E-06&	1.3621E-05&	&	6.0477E-06&	8.4966E-06&	2.1728E-05	\\
5	&1.1562E-03	&2.9942E-03	&1.2767E-03		&&4.4231E-06	&5.3967E-06&	5.3723E-05&	&	1.4659E-06&	2.2494E-06&	1.7805E-05	\\
7	&7.2867E-04	&1.8781E-03	&3.1572E-03		&&3.2083E-06	&3.7144E-06&	1.5291E-05&	&   4.7870E-06&	6.6079E-06&	2.2815E-05	\\
10	&5.8889E-04	&1.5158E-03	&2.2874E-03		&&3.1815E-06	&3.7741E-06&	1.3593E-05&	&	5.2716E-06&	7.3166E-06&	2.2524E-05\\
\bottomrule
\end{tabular}
\end{table}

 \begin{table}[!htbp]
\caption{Comparison of the error norms for Problem \ref{ex1} for $h=0.05,\triangle t=0.001$ }\label{tab-ex1.2}
\vspace{.1cm}
\centering
\begin{tabular}{*8l*4l}
\toprule
{}& \multicolumn{3}{l}{MCB-DQM (SSP-RK43) \cite{MB14} }  &{}& \multicolumn{7}{l}{MTB-DQM}   \\ \cline{2-4} \cline{6-12}
{}&\multicolumn{3}{c}{   }  &{}& \multicolumn{3}{l}{SSP-RK54 } & & \multicolumn{3}{l}{SSP-RK43} \\ \cline{2-4} \cline{6-8} \cline{10-12}
{t} &$L_2$ &$L_\infty$& $R_e$ &{}&$L_2$ &$L_\infty$ &$R_e$& &$L_2$ &$L_\infty$ &$R_e$\\
\midrule			
1	&9.8870E-05$\quad$ &	2.4964E-04$\quad$&	6.2977E-04&&	3.5763E-07$\quad$&	5.8782E-07$\quad$&	2.3919E-06&&	6.1937E-07$\quad$&	 7.4249E-07$\quad$&	 4.1424E-06\\
2	&1.2148E-04&	3.2296E-04&	1.0025E-03&&	4.4998E-07&	6.7260E-07&	3.8988E-06&&	3.4717E-07&	5.9880E-07&	3.0080E-06\\
3	&3.7627E-05&	9.9310E-05&	1.3078E-04&&	7.8209E-07&	1.2239E-06&	2.8543E-06&&	9.2093E-07&	1.2919E-06&	3.3620E-06\\
5	&1.2762E-04&	3.3205E-04&	1.5411E-03&&	3.6762E-07&	4.4816E-07&	4.6831E-06&&	2.1887E-07&	3.4820E-07&	2.7882E-06\\
7	&6.7672E-05&	1.7679E-04&	3.0892E-04&&	5.0460E-07&	8.8119E-07&	2.4186E-06&&	7.3202E-07&	1.0139E-06&	3.5087E-06\\
10	&5.1764E-05&	1.3521E-04&	2.1245E-04&&	5.8060E-07&	9.9247E-07&	2.5020E-06&&	8.0485E-07&	1.1198E-06&	3.4685E-06\\
\bottomrule
\end{tabular}
\end{table}

\begin{table}[!htbp]
\caption{Comparison of the error norms for Problem \ref{ex2} with $\alpha=10,\beta=5$ and $\Delta t=0.01$ and $h=0.1$}\label{tab-ex2.1}
\centering
\begin{tabular}{*9l*3l}
\toprule
$t$&  \multicolumn{7}{l}{MTB-DQM}  &{}& \multicolumn{3}{l}{MCB-DQM \cite{MB14}}   \\ \cline{2-8} \cline{10-12}
&  \multicolumn{3}{l}{SSP-RK54} & {}& \multicolumn{3}{l}{SSP-RK43}  &{}& \multicolumn{3}{l}{SSP-RK43}   \\ \cline{2-4}\cline{6-8} \cline{10-12}
{}&$L_2$ &$L_\infty$ & $R_e$ & {} &$L_2$ &$L_\infty$ & $R_e$&{}&$L_2$ &$L_\infty$ & $R_e$\\
\midrule
0.5~~&	6.9121E-06$~~$&	1.0190E-05$~~$&	2.6182E-05$~~$&{}&	4.4634E-06$~~$&	7.4933E-06$~~$&	1.6906E-05&&	8.3931E-04$~~$&	3.3019E-03$~~$&	2.8902E-03\\
1&	    5.3610E-06&	7.0248E-06&	3.3480E-05&&	3.1263E-06&	4.9241E-06&	1.9523E-05&&	6.0254E-04&	2.0597E-03&	3.4208E-03\\
2&	    2.2372E-06&	2.8580E-06&	3.7979E-05&&	1.2315E-06&	1.8734E-06&	2.0905E-05&&	2.4167E-04&	7.6531E-04&	3.7297E-03\\
3&	    8.3506E-07&	1.0602E-06&	3.8921E-05&&	4.5449E-07&	6.8641E-07&	2.1183E-05&&	8.9534E-05&	2.7920E-04&	3.7937E-03\\
5&	    1.1369E-07&	1.4412E-07&	3.9156E-05&&	6.1709E-08&	9.3030E-08&	2.1252E-05&&	1.2168E-05&	3.7800E-05&	3.8097E-03\\
\bottomrule
\end{tabular}
\end{table}

 \begin{table}[!htbp]
\caption{Comparison of the error norms for Problem \ref{ex2} with $\delta t=0.001,\alpha=10, \beta=0, 5$ and $h=0.05$}\label{tab-ex2.2}
\vspace{.1cm}
\centering
 \selectfont{\fontsize{8.0}{09.00}
\begin{tabular}{*9c*5c}
\toprule
$t$&  \multicolumn{7}{l}{ MTB-DQM }  &{}& \multicolumn{3}{l}{ MCB-DQM \cite{MB14} } &{}& \multicolumn{1}{l}{PDQM \cite{JPM12} }  \\ \cline{2-8} \cline{10-12} \cline{14-14}
{}&  \multicolumn{3}{l}{ SSP-RK54 } & {} & \multicolumn{3}{l}{ SSP-RK43 }  &{}& \multicolumn{3}{l}{SSP-RK43} &{}& \multicolumn{1}{l}{ }  \\ \cline{2-4}\cline{6-8}  \cline{10-12} \cline{14-14}
$\beta=5$ &$L_2$&	$L_\infty$&	$R_e$ & {} &$L_2$&	$L_\infty$&	$R_e$ & &	$L_2$&	$L_\infty$	 &$R_e$&    &	$R_e$ \\
\midrule																																				
0.5		&8.1325E-07$~$	&1.3162E-06	$~$&3.1852E-06&&	6.9072E-07$~~$&	1.2225E-06$~$&	2.7052E-06&&	1.0690E-04$~$&	2.4738E-04$~~$&	1.1088E-04&&	 1.1185E-04\\
1		&5.8464E-07	&8.4033E-07	&3.7753E-06&&	4.7197E-07&	7.6876E-07&	3.0477E-06&&	1.5293E-05&	3.3082E-04&	1.3266E-04&&	1.8051E-04\\
2		&2.3521E-07	&3.2221E-07	&4.1327E-06&&	1.8422E-07&	2.8699E-07&	3.2369E-06&&	4.6468E-05&	1.1380E-05&	3.1954E-04&&	4.7289E-04\\
3		&8.8084E-08	&1.1945E-07	&4.2070E-06&&	6.8575E-08&	1.0587E-07&	3.2753E-06&&	2.1994E-05&	4.3577E-05&	1.3024E-04&&	1.2656E-04\\
5		&1.1986E-08	&1.6213E-08	&4.2260E-06&&	9.3176E-09&	1.4353E-08&	3.2850E-06&&	2.7151E-06&	5.4141E-06&	1.4439E-04&&	9.2770E-04\\
$\beta=0$& &&&&&&&&&&& &  \\											
0.5		&7.2901E-07	&9.0802E-07	&2.8553E-06&&	8.2389E-07&	1.3548E-06&	3.2269E-06&&	9.2959E-05&	4.2348E-04&	3.4675E-04&&	1.1198E-04\\
1		&8.0718E-07	&1.0267E-06	&5.2123E-06&&	6.6905E-07&	9.4622E-07&	4.3200E-04&&	6.3652E-05&	2.5838E-04&	3.9146E-04&&	1.8635E-04\\
2		&5.7505E-07	&7.2597E-07	&1.0104E-05&&   3.5621E-07&	4.4306E-07&	6.2589E-06&&	2.5540E-05&	9.5843E-05&	4.2739E-04&&	5.1797E-04\\
3		&3.1143E-07	&3.9325E-07	&1.4874E-05&&	1.7331E-07&	2.1105E-07&	8.2775E-06&&	9.9234E-06&	3.5340E-05&	4.5140E-04&&	1.4412E-04\\
5		&6.7774E-08	&8.5735E-08	&2.3894E-05&&	3.5157E-08&	4.3088E-08&	1.2395E-05&&	1.5116E-06&	4.8043E-06&	5.0758E-04&&	1.0883E-04\\
\bottomrule
\end{tabular} }
\end{table}

 \begin{table}[!htbp]
\caption{Comparison of the error norms for Problem \ref{ex4} with $\alpha=10, 50,\beta=5$ and $h=0.05, \triangle t =0.001$}
\label{tab-ex4.1}
\centering
 \selectfont{\fontsize{8.0}{09.00}
\begin{tabular}{*9c*5c}
\toprule
$t$&  \multicolumn{7}{l}{ MTB-DQM }  &{}& \multicolumn{3}{l}{ MCB-DQM \cite{MB14} } &{}& \multicolumn{1}{l}{PDQM \cite{JPM12} }  \\ \cline{2-8} \cline{10-12} \cline{14-14}
{}&  \multicolumn{3}{l}{ SSP-RK43 } & {} & \multicolumn{3}{l}{ SSP-RK54 }  &{}& \multicolumn{3}{l}{SSP-RK43} &{}& \multicolumn{1}{l}{ }  \\ \cline{2-4}\cline{6-8}  \cline{10-12} \cline{14-14}
$\alpha=10$ &$L_2$ &	$L_\infty$ &$R_e$ & &$L_2$ &	$L_\infty$ &$R_e$ & &$L_2$ &	$L_\infty$ &$R_e$ & & $R_e$ \\
\midrule																			
0.5	&	9.2007E-07$~~~~$&	1.6964E-06$~~~~$&	2.4906E-06&&		2.0867E-06$~~~~$&	2.8539E-06$~~~~$&	5.6486E-06&&		1.0696E-04$~~~~$&	 3.7559E-04$~~~~$&	1.3787E-05&&		2.0149E-05\\
1	&	7.0486E-07&	1.1533E-06&	3.0991E-06&&		2.5049E-06&	3.2486E-06&	1.1014E-05&&		1.7174E-04&	5.6395E-04&	3.5957E-05&&		4.6003E-05\\
2	&	2.8961E-07&	6.3703E-07&	1.6496E-06&&		1.3896E-06&	1.7942E-06&	7.9150E-06&&		1.6468E-04&	5.1298E-04&	4.4722E-05&&		4.4460E-05\\
3	&	9.9370E-07&	1.8479E-06&	2.3841E-06&&		1.4014E-06&	2.2266E-06&	3.3622E-06&&		8.9858E-06&	1.9564E-05&	1.0266E-06&&		6.4830E-06\\
5	&	1.7624E-07&	3.7604E-07&	1.4759E-06&&		1.6566E-06&	2.1478E-06&	1.3873E-05&&		1.7737E-04&	5.5627E-04&	7.0735E-05&&		7.3626E-05\\
7	&	9.0389E-07&	1.5345E-06&	2.8482E-06&&		2.5349E-06&	3.2890E-06&	7.9874E-06&&		1.4200E-04&	4.7231E-04&	2.1307E-05&&		\\
10	&	1.4394E-07&	1.9332E-07&	2.3501E-06&&		1.8984E-06&	2.4642E-06&	6.9185E-06&&		1.2241E-04&	4.1222E-04&	1.6514E-05&&		2.4901E-05\\
$\alpha=50$	& &&&&&&&&& \\													
0.5	&	5.9234E-07&	1.3746E-06&	1.6034E-06&&		2.2131E-06&	3.2840E-06&	5.9907E-06&&		9.8800E-05&	3.6962E-04&	1.2735E-05&&		2.8724E-05\\
1	&	5.9130E-07&	1.1709E-06&	2.5998E-06&&		3.2436E-06&	4.4896E-06&	1.4261E-05&&		1.6766E-04&	5.6874E-04&	3.5104E-05&&		7.0064E-05\\
2	&	2.0351E-07&	3.0955E-07&	1.1592E-06&&		2.4070E-06&	3.2577E-06&	1.3710E-05&&		1.7109E-04&	5.2572E-04&	4.6408E-05&&		6.7759E-05\\
3	&	6.6949E-07&	1.5167E-06&	1.6063E-06&&		1.6272E-06&	2.8374E-06&	3.9041E-06&&		1.7406E-05&	4.3459E-05&	1.9886E-06&&		1.3856E-05\\
5	&	2.9450E-07&	4.5993E-07&	2.4663E-06&&		3.1533E-06&	4.2936E-06&	2.6408E-05&&		1.8420E-04&	5.6940E-04&	7.3460E-05&&		1.2840E-04\\
7	&	7.7456E-07&	1.5214E-06&	2.4407E-06&&		3.5804E-06&	5.1320E-06&	1.1282E-05&&		1.3760E-04&	4.7587E-04&	2.0647E-05&&		\\
10	&	7.9206E-07&	1.6017E-06&	2.2439E-06&&		3.3625E-06&	4.8643E-06&	9.5260E-06&&		1.1691E- 4&	4.1396E-04&	1.5772E-05&&		4.4202E-05\\
\bottomrule
\end{tabular} }
\end{table}

\begin{table}[!htbp]
\caption{Comparison of the error norms for Problem \ref{ex5} with $\alpha=1,\beta=1$ for different values of $\triangle t, h$}\label{tab5.1}
\vspace{.2cm}
\centering
\begin{tabular}{*5l*7l}
\toprule
$t$&  \multicolumn{7}{l}{ MTB-DQM: $\triangle t =0.1h=0.01$ }  &{}& \multicolumn{3}{l}{ MCB-DQM: $\triangle t =0.1h=0.01$ \cite{MB14}}   \\ \cline{2-8} \cline{10-12}
{}&  \multicolumn{3}{l}{ SSP-RK54}  &{}& \multicolumn{3}{l}{ SSP-RK43} & {}& \multicolumn{3}{l}{ SSP-RK43}   \\ \cline{2-4} \cline{6-8}  \cline{10-12}
{}	&$L_2$ &$L_\infty$ & $R_e$ &{}&$L_2$ &$L_\infty$ &$R_e$ &{}&$L_2$ &$L_\infty$ &$R_e$  \\
\midrule
1&	3.9792E-04$\quad$&	6.7060E-04$\quad$&	3.2826E-04&&	4.5279E-04$\quad$&	7.3999E-04$\quad$&	3.7353E-04&&	1.4441E-02$\quad$&	2.9996E-02$\quad$&	 1.0829E-03\\
2&	4.5063E-05&	1.1081E-04&	1.0105E-04&&	4.7152E-05&	1.1332E-04&	1.0574E-04&&	1.3898E-03&	3.9711E-03&	2.8333E-04\\
3&	4.0577E-05&	7.4504E-05&	2.4983E-04&&	4.5408E-05&	8.0118E-05&	2.7957E-04&&	1.3018E-03&	2.2178E-03&	7.2867E-04\\
5&	4.1054E-06&	8.6394E-06&	1.8677E-04&&	4.3519E-06&	8.9406E-06&	1.9798E-04&&	1.1112E-04&	2.0618E-04&	4.5956E-04\\
7&	4.6720E-07&	1.0446E-06&	1.5705E-04&&	4.7606E-07&	1.0625E-06&	1.6003E-04&&	1.3695E-05&	3.0052E-05&	4.1851E-04\\
10&	3.8677E-08&	7.0406E-08&	2.6114E-04&&	4.3364E-08&	7.5959E-08&	2.9278E-04&&	1.4408E-06&	2.5354E-06&	8.8440E-04\\
\multicolumn{12}{l}{ $h=.05, \triangle t=0.001$}\\								
0.5&1.2789E-04&	2.6720E-04&	6.4982E-05&&	1.2998E-04&	2.6974E-04&	6.6047E-05&&	3.4808E-03&	9.5129E-03&	8.4225E-05\\
1&	1.0533E-04&	1.8245E-04&	8.8240E-05&&	1.0829E-04&	1.8630E-04&	9.0721E-05&&	3.2351E-03&	7.4749E-03&	1.2906E-04\\
2&	1.0444E-05&	3.0713E-05&	2.3809E-05&&	1.0524E-05&	3.0837E-05&	2.3990E-05&&	2.8518E-04&	1.0361E-03&	3.0957E-05\\
3&	1.0890E-05&	2.0527E-05&	6.7481E-05&&	1.1155E-05&	2.0822E-05&	6.9120E-05&&	3.1028E-04&	5.7859E-04&	9.1555E-05\\
5&	1.0291E-06&	2.3991E-06&	4.7070E-05&&	1.0364E-06&	2.4103E-06&	4.7404E-05&&	2.4495E-05&	6.7234E-05&	5.3354E-05\\
7&	1.0346E-07&	2.5898E-07&	3.4969E-05&&	1.0290E-07&	2.5876E-07&	3.4780E-05&&	2.5376E-06&	8.2203E-06&	4.0840E-05\\
10&	1.1780E-08&	2.1778E-08&	8.0050E-05&&	1.2098E-08&	2.2179E-08&	9.0721E-05&&	3.6505E-06&	8.5897E-06&	1.1812E-05\\
\bottomrule
\end{tabular}
\end{table}

\begin{table}[!htbp]
\caption{Comparison of the error norms for Problem \ref{ex6} with $\alpha=\beta=1$ for different values of $\triangle t, h$}\label{tab-ex6.1}
\centering
\begin{tabular}{*5l*4l}
\toprule
$t$&  \multicolumn{5}{l}{ MTB-DQM: $\triangle t =0.1h=0.01$  }  &{}& \multicolumn{2}{l}{ MCB-DQM: $\triangle t =0.1h=0.01$ \cite{MB14}}   \\ \cline{2-6} \cline{8-9}
{}&  \multicolumn{2}{l}{ SSP-RK54}  &{}& \multicolumn{2}{l}{ SSP-RK43} & {}& \multicolumn{2}{l}{ SSP-RK43}   \\ \cline{2-3} \cline{5-6}  \cline{8-9}
{}	&$L_2$ &$L_\infty$ &{}&$L_2$ &$L_\infty$ &{}&$L_2$ &$L_\infty$\\
\midrule									
1$\quad$&  5.7445E-04$\quad$$\qquad$&	7.1729E-04$\qquad$&&	3.5177E-04$\quad$$\qquad$&	4.5346E-04$\qquad$&	&1.6144E-03$\quad$$\qquad$&	3.6006E-03\\
2&  1.7401E-04&	2.2454E-04&&	1.0853E-04&	1.4837E-04&&	2.6345E-03&	5.7068E-03\\
3&  1.9329E-05&	2.1560E-05&&	1.2187E-05&	1.5671E-05&&	5.3845E-04&	1.2479E-03\\
5&	6.5016E-06&	8.5888E-06&&	4.0448E-06&	5.6143E-06&&	1.2418E-04&	2.1003E-04\\
7&	1.3046E-06&	1.6306E-06&&	8.0131E-07&	1.0392E-06&&	1.3653E-05&	2.6261E-05\\
10&	5.8352E-08&	7.3725E-08&&	3.5842E-08&	4.6819E-08&&	7.5592E-06&	1.4083E-06\\
\multicolumn{9}{l}{ $h=.05, \triangle t=0.001$}\\	
0.5&3.2619E-05&	4.6318E-05&&	2.6638E-05&	4.0745E-05&&	3.5833E-04&	9.5129E-04\\
1&	5.5101E-05&	7.2244E-05&&	4.3591E-05&	5.8390E-05&&	3.2351E-04&	7.4749E-04\\
2&	1.5540E-05&	2.1383E-05&&	1.2429E-05&	1.7899E-05&&	2.8518E-05&	1.0361E-04\\
3&	8.3614E-07&	1.1513E-06&&	7.4268E-07&	1.1493E-06&&	3.1028E-05&	5.7859E-04\\
5&	5.1816E-07&	7.4785E-07&&	4.1500E-07&	6.2203E-07&&	2.4495E-06&	6.7234E-05\\
7&	1.5582E-07&	1.9985E-07&&	1.2341E-07&	1.6197E-07&&	2.5376E-07&	8.2203E-07\\
10&	7.1283E-09&	9.0996E-09&&	5.6224E-09&	7.2378E-09&&	3.6505E-09&	8.5897E-08\\
\bottomrule
\end{tabular}
\end{table}

\begin{table}[!htbp]
\caption{Comparison of the error norms for Problem \ref{ex7} with $\alpha=\beta=1$ taking $h=0.05, \triangle t=0.001$}\label{tab-ex7.0}
\vspace{.2cm}
\centering
\begin{tabular}{*9c*6c}
\toprule
$t$&  \multicolumn{3}{l}{MTB-DQM (SSP-RK43) }  &{}& \multicolumn{3}{l}{MCB-DQM (SSP-RK43)\cite{MB14}} &{}& \multicolumn{2}{l}{ Dehghan and Ghesmati \cite{DG10a}}  \\ \cline{2-4} \cline{6-8} \cline{10-11}
{} &	$L_2$ & $L_\infty$ &	$R_e$ &&	$L_2$ &	$L_\infty$ &	 $R_e$ &&	$R_e:$MLWS &	$R_e:$MLPG\\
\midrule					
0.5&	4.8136E-05$\quad$&	9.7584E-05$\quad$&	5.2424E-05&&	1.0690E-03$\quad$&	2.4738E-03$\quad$&	1.1088E-03&&	7.9390E-05&	9.9910E-05\\
1&	7.3094E-05&	1.0833E-04&	6.6578E-05&&	1.5293E-03&	3.3082E-03&	1.3266E-03&&	9.0980E-05&	7.1980E-05\\
2&	2.9279E-05&	4.9102E-05&	2.1140E-05&&	4.6468E-04&	1.1380E-03&	3.1954E-04&&	8.7050E-04&	8.7840E-05\\
3&	1.1806E-05&	2.0621E-05&	7.3411E-06&&	2.1994E-04&	4.3577E-04&	1.3024E-04&&	9.9310E-04&	4.8010E-04\\
4&	1.2638E-05&	1.9269E-05&	7.0574E-06&&	2.7151E-04&	5.4141E-04&	1.4435E-05&&	4.7030E-03&	6.0910E-04\\
5&	7.8175E-06&	1.2251E-05&	4.0193E-06&&	1.7201E-04&	3.4812E-04&	8.4225E-05&&	7.3020E-03&	9.4980E-04\\
7&	4.8159E-06&	7.4804E-06&	2.1924E-06&&	\multicolumn{6}{l}{} \\				
10&	2.6185E-06&	4.0492E-06&	1.0539E-06&&	7.7285E-05&	1.4037E-04&	2.9624E-05&\multicolumn{3}{l}{}\\	
\bottomrule
\end{tabular}
\end{table}

\end{document}